\documentclass[leqno,a4paper,12pt]{elsart}

\usepackage{amssymb}
\usepackage{latexsym}
\usepackage{amsmath}

\newtheorem{theorem}{Theorem}[section]
\newtheorem{proposition}[theorem]{Proposition}
\newtheorem{corollary}[theorem]{Corollary}
\newtheorem{lemma}[theorem]{Lemma}
\newtheorem{defi2}[theorem]{Definition}
\newtheorem{remarke}[theorem]{Remark}
\newtheorem{assump}[theorem]{Assumption}
\newtheorem{exx}[theorem]{Example}

\newenvironment{remark}{\begin{remarke}\rm}{\end{remarke}}

\newcommand{\R}{{\mathbb R}}
\newcommand{\natt}{{\mathbb N}}

\newcommand{\norm}[1]{\left\| #1 \right\|}

\newcommand{\ind}{\mbox{{\Large $\chi$}}}

\newcommand{\BMO}{\mbox{\rm BMO}}

\newcommand{\Comp}{{\mathbb C}}

\newcommand{\g}{{\mathcal G}}
\newcommand{\area}{{\mathcal A}}
\newcommand{\Tt}{{\mathcal T}}
\newcommand{\Pt}{{\mathcal P}}

\newcommand{\pa}{\partial}
\newcommand{\F}{{\mathcal F}}

\newcommand{\B}{{\mathbb B}}
\newcommand{\T}{{\mathbb T}}
\newcommand{\AB}{{\mathbb A}}
\newcommand{\D}{{\mathbb D}}

\newcommand{\cqd}{\hfill$\Box$}

\numberwithin{equation}{section}

\begin{document}


\begin{frontmatter}



\title{Vector-valued Littlewood-Paley-Stein theory for semigroups}

\author{Teresa Mart\'\i nez\thanksref{madrid}}
\address{Departamento de Matem\'aticas, Facultad de Ciencias, Universidad Au\-t\'o\-noma de
Madrid, 28049 Madrid, Spain}
\ead{teresa.martinez@uam.es}

\author{Jos\'e L. Torrea\thanksref{madrid}}
\address{Departamento de Matem\'aticas, Facultad de Ciencias, Universidad Au\-t\'o\-noma de
Madrid, 28049 Madrid, Spain}
\ead{joseluis.torrea@uam.es}

\author{Quanhua Xu}
\address{Universit{\'e} de Franche-Comt{\'e}, Equipe de Math{\'e}matiques,
UFR des Sciences et Techniques, 16, Route de Greay,
F-25030 Besan\c{c}on Cedex, France}
\ead{qx@math.univ-fcomte.fr}

\thanks[madrid]{First and second authors were partially supported by
RTN Harmonic Analysis and Related Problems contract HPRN-CT-2001-00273-HARP
and by BFM grant 2002-04013-C02-02}


\begin{abstract}
We develop a generalized Littlewood-Paley theory for semigroups
acting on $L^p$-spaces of functions with values in uniformly
convex or smooth Banach spaces. We characterize, in the
vector-valued setting, the validity of the one-sided inequalities
concerning the generalized Littlewood-Paley-Stein $g$-function
associated with a subordinated Poisson symmetric diffusion
semigroup by the martingale cotype and type properties of the
underlying Banach space. We show that in the case of the usual
Poisson semigroup and the Poisson semigroup subordinated to the
Ornstein-Uhlenbeck semigroup on $\R^n$, this general theory
becomes more satisfactory (and easier to be handled) in virtue of
the theory of vector-valued Calder{\'o}n-Zygmund singular integral
operators.
\end{abstract}

\begin{keyword}
Littlewood-Paley theory, semigroups, uniformly convex or smooth
Banach spaces, vector-valued Calder{\'o}n-Zygmund operators.

\smallskip
\MSC{46B20; 42B25,42A61.}

\end{keyword}

\end{frontmatter}


\section{Introduction and preliminaries}

Given a martingale $\{f_n\}$ with values in a Banach space $\B$, its generalized
\lq\lq square\rq\rq\ function is defined as
$$
S_qf=\Big(\sum_{n=1}^\infty\|f_n-f_{n-1}\|^q_\B\Big)^{1/q}.
$$
Then  $\B$ is said to have martingale cotype $q$, $2\leq q<\infty$
if there exist $p\in(1,\infty)$ and a constant $C>0$ such that
$\|S_qf\|_{L^p}\leq C\sup_n\|f_n\|_{L^p_\B}$ for every bounded
$\B$-valued $L^p$-martingale $\{f_n\}$. The validity of the
reverse inequality defines martingale type $q$, $1< q\leq 2$.
Recall that if the inequality above (or its inverse) holds for one
$p\in(1,\infty)$, so does it for all $p\in(1,\infty)$. These
notions were introduced and studied in depth by Pisier (see
\cite{pi:75,pi:86}). His renorming theorem states that there are
geometric properties of the underlying Banach space, characterized
by the existence of an equivalent norm in the space which is
uniformly convex of power type $q$ or uniformly smooth of power
type $q$. We also recall that $\B$ is of martingale cotype $q$ iff
$\B^*$ is of martingale type $q'$, where $q'$ is the index
conjugate to $q$. For $p\in(1,\infty)$, $L^p$ is of martingale
cotype $\max\{2,p\}$ and of martingale type $\min\{2,p\}$.

\medskip

On the other hand, it is well known that martingale inequalities
involving square function are closely related to the corresponding
inequalities concerning the Li\-ttle\-wood-Paley or Lusin square
function in harmonic analysis. It is in this spirit that a
generalized Littlewood-Paley theory is developed in \cite{xu} for
functions with values in uniformly convex Banach spaces. Let us
recall the main results of \cite{xu}.  Let $f$ be a function in
$L^1(\T)$, where $\T$ denotes the torus equipped with normalized
Haar measure $d\theta$. The classical Littlewood-Paley
$g$-function is defined for $z\in\T$ as
\begin{eqnarray*}
Gf(z)&=&\Big(\int_0^1 (1-r)^2\|\nabla
P_r*f(z)\|^2\,\frac{dr}{1-r}\Big)^{1/2}.
\end{eqnarray*}
In this notation,
\begin{equation}
\label{gradient T} \|\nabla P_r*f(z)\| =
\Big(\big|\frac{\partial P_r}{\partial r}*f(z)\big|^2+
\big|\frac{1}{r}\frac{\partial P_r}{\partial
\theta}*f(z)\big|^2\Big)^{1/2},
\end{equation}
with
$$
P_r(\theta)=\frac{1-r^2}{1+r^2-2r\cos \theta},
$$
being the Poisson kernel for the disk. It is a classical result
that for any $p\in(1,\infty)$ there exist positive constants $c_p$
and $C_p$ such that
\begin{equation}\label{equivalence}
c_p\|f\|_{L^p(\T)}\leq|\hat{f}(0)|+\|Gf\|_{L^p(\T)}\leq
C_p\|f\|_{L^p(\T)}.
\end{equation}
One can extend the definition of $G$ to functions that take values
in a Banach space $\B$, just by replacing absolute value by norm
in (\ref{gradient T}). In this case, (\ref{equivalence}) holds if
and only if $\B$ is isomorphic to a Hilbert space. However, one of
the two inequalities in (\ref{equivalence}) can be true in non
Hilbertian spaces. The study of these one-sided inequalities is
the main objective of \cite{xu}.  More generally, we can introduce
the following generalized \lq\lq Littlewood-Paley
$g$-function\rq\rq\
$$
G_qf(z)=\Big(\int_0^1 (1-r)^q\|\nabla
P_r*f(z)\|_{\B}^q\,\frac{dr}{1-r}\Big)^{1/q}\ .
$$
Then $\B$ is said to be of Lusin cotype $q$ (resp. Lusin type $q$)
if there exist $p\in(1,\infty)$ and a positive constant $C$ such
that
$$
\|G_qf\|_{L^p(\T)}\leq C\|f\|_{L^p_\B(\T)}\quad \Big({\rm resp.}\
\|f\|_{L^p_\B(\T)}\leq
C\big(\|\hat{f}(0)\|_\B+\|G_qf\|_{L^p(\T)}\big)\Big) .
$$
It is not difficult to see that if $\B$ is of Lusin cotype $q$
(resp. Lusin type $q$), then $2\leq q\leq\infty$ (resp. $1\leq
q\leq 2$). It is proved in \cite{xu} that the definition above is
independent of $p$, that is,  if one of the inequalities above
holds for one $p\in(1,\infty)$, then so does it for every
$p\in(1,\infty)$ (with a different constant depending on $p$). The
main result of \cite{xu} states that a Banach space $\B$ is of
Lusin type $q$ (resp. Lusin cotype $q$) iff $\B$ is of martingale
type $q$ (resp. martingale cotype $q$).

\medskip

The main goal of the present paper is to extend the results in
\cite {xu} to general symmetric diffusion semigroups, and thus to
develop a generalized Littlewood-Paley theory for these semigroups
on $L^p$-spaces of functions with values in uniformly convex or
smooth Banach spaces. Recall that a symmetric diffusion semigroup
is a collection of linear operators $\{\Tt_t\}_{t\geq 0}$ defined
on $L^p(\Omega,d\mu)$ over a measure space $(\Omega,d\mu)$
satisfying the following properties
\begin{equation}\label{semigroup 1}
\Tt_0=\mbox{ Id},\quad \Tt_t\Tt_s=\Tt_{t+s},\quad
\|\Tt_t\|_{L^p\to L^p}\leq 1, \ \forall\;p\in[1,\infty];
\end{equation}
\begin{equation}\label{semigroup 2}
\lim_{t\rightarrow 0}\Tt_tf=f\mbox{ in }L^2,\quad \forall\;f\in
L^2;
\end{equation}
\begin{equation}\label{semigroup 3}
\Tt_t^*=\Tt_t\mbox { on }L^2,\quad  \Tt_tf\geq 0\mbox{ if }f\geq
0, \quad \Tt_t1=1.
\end{equation}
The subordinated Poisson semigroup $\{\Pt_t\}_{t\geq 0}$ is
defined as
\begin{equation}\label{subordinated semigroup}
\Pt_tf=\frac{1}{\sqrt\pi}\int_0^\infty \frac{e^{-u}}{\sqrt
u}\Tt_{t^2/4u}f\,du= \frac{t}{2\sqrt\pi}\int_0^\infty
\frac{e^{-t^2/4u}}{u^{3/2}}\Tt_uf\,du.
\end{equation}
$\{\Pt_t\}_{t\geq 0}$ is again a symmetric diffusion semigroup,
see \cite{st:}. Recall that if $A$ denotes the infinitesimal
generator of $\{\Tt_t\}_{t\geq 0}$, then that of $\{\Pt_t\}_{t\geq
0}$ is $-(-A)^{1/2}$.

\medskip

It is well know (and easy to check) that any bounded operator $T$
on $L^p(\Omega)$ for all $p\in[1, \infty]$ naturally and boundedly
extends to $L^p_\B(\Omega)$ for every Banach space $\B$, where
$L^p_\B(\Omega)$ denotes the usual Bochner-Lebesgue $L^p$-space of
$\B$-valued functions defined on $\Omega$. More precisely, the
extension is $T\otimes\mbox{Id}_{\B}$. Indeed, this is clear for
$p=1$ (via projective tensor product); the case $p=\infty$ is done
by duality, and then $1<p<\infty$ by interpolation. With a slight
abuse of notation (which will not cause any ambiguity), we shall
denote these extensions still by the same symbol $T$.

Thus $\Tt_t$ and $\Pt_t$  have straightforward extensions to
$L^p_\B(\Omega)$ for every Banach space $\B$; moreover, these
extensions are also contractive. (Note that we can also justify
these extensions by the positivity of $\Tt_t$ and $\Pt_t$.)
According to the convention above, we shall consider
$\{\Tt_t\}_{t\geq 0}$ and $\{\Pt_t\}_{t\geq 0}$ as semigroups on
$L^p_\B(\Omega)$ too.

\medskip

In these circumstances we can define the generalized \lq\lq
Littlewood-Paley $g$-function\rq\rq\ associated to the semigroup
as
$$
{\mathfrak G}_q(f)(x)=\Big(\int_0^\infty \Big\|t\frac{\partial
\Pt_tf(x)}{\partial t} \Big\|_{\B}^q\,\frac{dt}{t}\Big)^{1/q}.
$$
The first result of this paper, see Theorem \ref{lusin
semigroups},  states that a Banach space $\B$ is of martingale
cotype  $q$ iff for every symmetric diffusion semigroup
$\{\Tt_t\}_{t\geq 0}$ with subordinated semigroup
$\{\Pt_t\}_{t\geq 0}$, the generalized $g$-function operator
${\mathfrak G}_q$ is bounded from $L^p_\B(\Omega)$ to
$L^p(\Omega)$, namely
\begin{equation}\label{cotype inequality semigroup}
 \|{\mathfrak G}_q(f)\|_{L^p(\Omega)}\le C\|f\|_{L_\B^p(\Omega)}\ ,
 \quad\forall\; f\in L_\B^p(\Omega).
\end{equation}
The validity of the reverse inequality (with a necessary
additional term) characterizes the martingale type $q$ (see
Theorem \ref{lusin type semigroups}). These results are proved in
section 2. The main ingredient of our arguments is the classical
Rota theorem on the dilation of a positive contraction on $L^p$ by
conditional expectations. This theorem allows to reduce
(\ref{cotype inequality semigroup}) (after a discretization) to a
corresponding inequality for martingales.

\medskip

This approach via Rota's theorem is also efficacious in studying
(\ref{cotype inequality semigroup}) and its dual form for an
individual semigroup. We shall show in section 3 that for a given
subordinated Poisson semigroup $\{\Pt_t\}$, (\ref{cotype
inequality semigroup}) is equivalent to its dual form, which is an
inequality reverse to (\ref{cotype inequality semigroup}) with
$\B, p$ and $q$ replaced by $\B^*, p'$ and $q'$, respectively (and
with an additional term). The key to this is the existence of a
certain projection, whose proof, using Rota's theorem once more,
is unfortunately rather technical and complicated.

\medskip

Our proof for the implication \lq\lq(\ref{cotype inequality
semigroup}) $\Rightarrow$ martingale cotype $q$\rq\rq\ uses the
Poisson semigroup on the torus modulo the results in \cite {xu}
quoted previously. (Note however that this Poisson semigroup on
the torus is a multiplicative semigroup on $(0,1)$.) Thus it would
be interesting to know the family of  semigroups $\{\Pt_t\}_{t\geq
0}$ for which the validity of (\ref{cotype inequality semigroup})
implies martingale cotype $q$. One of the aims of the remainder of
the paper (after section 3) is to show that this is indeed the
case for the Poisson semigroup on $\R^n$. Such a result is, of
course, conceivable after \cite{xu}. For any $q\geq 1$, the
$n$-dimensional generalized \lq\lq Littlewood-Paley
$g$-function\rq\rq\ is defined as
$$
\g_q(f)(x)=\Big(\int_0^\infty t^q\|\nabla P_t*
f(x)\|_{\ell^2_\B}^q\,\frac{dt}{t}\Big)^{1/q},
$$
where
$$
\|\nabla P_t*f(x)\|_{\ell^2_\B} = \Big(\Big\|\frac{\partial
P_t}{\partial t}*f(x)\Big\|^2_\B+ \sum_{k=1}^n\Big\|\frac{\partial
P_t}{\partial x_k}*f(x)\Big\|^2_\B\Big)^{1/2},
$$
with
$$
P_t(x)=\frac{\Gamma(\frac{n+1}{2})}{\pi^{\frac{n+1}{2}}}
\frac{t}{(|x|^2+t^2)^{\frac{n+1}{2}}}\ ,
$$
the kernel of the Poisson semigroup for the upper half space. Note
that we use the same symbol $P_t$ to denote the Poisson kernels
both on $\T$ and on $\R^n$. This should not have any confusion in
the concrete context. Then $\B$ is of martingale cotype $q$ (resp.
martingale type $q$) iff for some (equivalently every) $p\in(1,
\infty)$ there is a constant $C$ such that
$$
 \|\g_q(f)\|_{L^p(\R^n)}\le C\|f\|_{L_\B^p(\R^n)}\quad
 \Big(\mbox{resp.}\;
 \|f\|_{L_\B^p(\R^n)}\le C  \|\g_q(f)\|_{L^p(\R^n)}\Big)
 \ .
$$
This result, among some others on $\g_q$-function on $\R^n$, is
proved in sections 4 and 5 and achieved by viewing the operators $\g$
as vector-valued Calder{\'o}n-Zygmund operators. For these operators, under
suitable conditions, one can get the equivalence of the strong type $(p,p)$ and
the boundedness $BMO-BMO$ (see Theorem \ref{singular integral}).
As a consequence, we obtain the characterization of the Lusin cotype
in terms of $BMO$ boundedness of the $g$-functions
(Corollary \ref{littlewood-paley} and Theorems \ref{cotype in R} and
\ref{type in R}). These two sections extend most of the
results in \cite{xu} for $\T$ to $\R^n$.

\medskip

The previous results for the usual Poisson semigroup on $\R^n$ can
be extended to the Poisson semigroup subordinated to the
Ornstein-Uhlenbeck semigroup on $\R^n$. This is done in section 6
(see Theorems \ref{cotype Orstein-Uhlenbeck}
and \ref{type Orstein-Uhlenbeck}).

\medskip

The last section  contains a further characterization of Lusin
cotype property in terms of almost sure finiteness of the
generalized Littlewood-Paley $g$-functions (Theorems \ref{finiteness ae th}
and \ref{martingale th}).


\section{One-sided vector-valued Littlewood-Paley-Stein
 inequalities for semigroups}

We shall consider general symmetric diffusion semigroups, that is,
the collections of linear operators  $\{\Tt_t\}_{t\geq 0}$ defined
on $L^p(\Omega)$, satisfying (\ref{semigroup 1}) - (\ref{semigroup
3}). Given such a semigroup $\{\Tt_t\}_{t\geq 0}$ we consider its
subordinated semigroup $\{\Pt_t\}_{t\geq 0}$, defined as in
(\ref{subordinated semigroup}). Let ${\mathbb F}\subset
L^2(\Omega)$ be the subspace of the fix points of
$\{\Pt_t\}_{t\geq 0}$, i.e., the subspace of all $f$ such that
$\Pt_t(f)=f$ for all $t>0$. Let $F:L^2(\Omega)\to {\mathbb F}$ be
the orthogonal projection. It is clear that $F$ extends to a
contractive projection (still denoted by $F$) on $L^p(\Omega)$ for
every $1\le p\le\infty$ and that $F\big(L^p(\Omega)\big)$ is
exactly the fix point space of $\{\Pt_t\}_{t\geq 0}$ on
$L^p(\Omega)$. Moreover, for any Banach space $\B$, $F$ extends to
a contractive projection on $L_B^p(\Omega)$ for every $1\le
p\le\infty$ and that $F\big(L_\B^p(\Omega)\big)$ is again the fix
point space of $\{\Pt_t\}_{t\geq 0}$ considered as a semigroup on
$L_\B^p(\Omega)$. According to our convention, in the sequel, we
shall use the same symbol $F$ to denote any of these contractive
projections.

Recall that the generalized Littlewood-Paley $g$-function
associated with $\{\Pt_t\}_{t\geq 0}$ is defined by
$$
{\mathfrak G}_q(f)(x)=\Big(\int_0^\infty \Big\|t\frac{\partial
\Pt_tf(x)}{\partial t} \Big\|_{\B}^q\,\frac{dt}{t}\Big)^{1/q}\ .
$$
The main results of this section are the following two theorems.

\begin{theorem}\label{lusin semigroups}
Given a Banach space $\B$ and $2\leq q<\infty$, the following
statements are equivalent:
\begin{itemize}
\item[{\sf i)}] $\B$ is of martingale cotype $q$.
\item[{\sf ii)}] For every symmetric diffusion semigroup
$\{\Tt_t\}_{t\geq 0}$ with subordinated semigroup
$\{\Pt_t\}_{t\geq 0}$ and for every (or, equivalently, for some)
$p\in(1,\infty)$ there is a constant $C$ such that
\begin{equation}\label{cotype semigroup}
\|{\mathfrak G}_qf\|_{L^p(\Omega)}\leq C\;\|f\|_{L^p_\B(\Omega)},
\quad \forall\;f\in L^p_\B(\Omega).
\end{equation}
\end{itemize}
\end{theorem}

\begin{theorem}\label{lusin type semigroups}
Given a Banach space $\B$ and $1< q\le 2$, the following
statements are equivalent:
\begin{itemize}
\item[{\sf i)}] $\B$ is of martingale type $q$.
\item[{\sf ii)}] For every symmetric diffusion semigroup
$\{\Tt_t\}_{t\geq 0}$ with subordinated semigroup
$\{\Pt_t\}_{t\geq 0}$ and for every (or, equivalently, for some)
$p\in(1,\infty)$ there is a constant $C$ such that
\begin{equation}\label{type semigroup}
\|f\|_{L^p_\B(\Omega)}\le C\big(\|F(f)\|_{L^p_\B(\Omega)} +
\|{\mathfrak G}_qf\|_{L^p(\Omega)}\big), \quad \forall\;f\in
L^p_\B(\Omega).
\end{equation}
\end{itemize}
\end{theorem}

The rest of this section is essentially devoted to the proof of
these theorems. The difficult part is the implication \lq\lq {\sf
i)} $\Rightarrow$ {\sf ii)}\rq\rq\ in Theorem \ref{lusin
semigroups}. Then the same implication in Theorem \ref{lusin type
semigroups} will follow by duality. Both converse implications
will be done by using the Poisson semigroup on the torus with the
help of \cite{xu}. For the main part of the proof we shall need
the following result, which has independent interest.

\begin{theorem}\label{theorem on the means}
Let $\B$ be a Banach space of martingale cotype $q\in[2,\infty)$
and $\{\Tt_t\}_{t\geq 0}$ a symmetric diffusion semigroup. Then
for any $p\in (1,\infty)$,
$$
\Big\|\Big(\int_0^\infty \Big\|t\frac{\partial M_tf}{\partial
t}\Big\|^q_\B \,\frac{dt}{t}\Big)^{1/q}\Big\|_{L^p(\Omega)} \leq
C_{p,q, \B}\; \|f\|_{L^p_\B(\Omega)}, \quad \forall\;f\in
L^p_\B(\Omega),
$$
where $\displaystyle M_t=\frac{1}{t}\int_0^t \Tt_s\,ds$.
\end{theorem}

The pattern of our proof for the theorem above is borrowed from
\cite [Chapter IV] {st:}. As in \cite {st:}, the key ingredient is
Rota's dilation theorem (see Theorem \ref{rota} below), which
allows to reduce the inequality in Theorem \ref{theorem on the
means} to a similar inequality for martingales.

Given a $\sigma$-finite measure space $(M,\F,dm)$ and a
sub-$\sigma$-algebra $\g\subset\F$,  we denote as usual by
$E(\,\cdot\,|\g)$ the conditional expectation with respect to
$\g$. (Note that our measure space $(M,\F,dm)$ is no longer a
probability one; however all usual properties on conditional
expectations in the probabilistic case are still valid in the
present setting.) Recall that $E(\,\cdot\,|\g)$ is a positive
contraction on $L^p(M,\F,dm)$ for every $p\in [1,\infty]$ and
naturally extends to $L_\B^p(M,\F,dm)$ for every Banach space
$\B$. The classical Doob maximal inequality is also valid in the
vector-valued setting. Let $(\F_n)$ be an increasing filtration of
sub-$\sigma$-algebras  of $\F$. For $f\in L_\B^1(M,\F,dm)$ we
define its maximal function as
$$
f^*=\sup_{n\geq 1}\|E(f|\F_n)\|_{\B}.
$$
Then we have the following  Doob maximal weak type (1,1)
inequality
$$
\lambda m\{f^*>\lambda\}\leq
\int_{\{f^*>\lambda\}}\|f(x)\|_\B\,dm(x)
$$
for every Banach space $\B$. Similarly,  we can also extend the
results of \cite{ma-to:00}; in particular, we get that for every
$1<p,q<\infty$ and every sequence $(f_k)\subset L_\B^p(M,\F,dm)$
\begin{equation}\label{p inequality}
\Big\|\Big(\sum_{k=1}^\infty((f_k)^*)^q\Big)^{1/q}\Big\|_{L^p}
\leq C_{p,q}
\Big\|\Big(\sum_{k=1}^\infty\|f_k\|_\B^q\Big)^{1/q}\Big\|_{L^p}.
\end{equation}
We shall use the following lemma, motivated by  \cite{st:}, inequality
$(**)$ on p.115.

\begin{lemma}\label{cesaro}
Let $\B$ be a Banach space of martingale cotype $q\in[2,\infty)$,
$(M,dm)$ be any $\sigma$-finite measure space and $\{E_n\}$ be an
arbitrary monotone sequence of conditional expectations on
$(M,dm)$. Then, for every $p$, $1<p<\infty$,
$$
\Big\|\Big(\sum_{n=1}^\infty n^{q-1}\|(\sigma_n-\sigma_{n-1})
f\|_\B^q\Big)^{1/q}\Big\|_{L^p} \leq C_{p,q, \B}\; \|f\|_{L^p_\B},
$$
where
$$
\quad\sigma_n=\frac{E_0+\cdots+E_n}{n+1}\ .
$$
\end{lemma}

\noindent{\sc Proof}.
Observe that it is enough to prove the inequality taking the summation
in $n\geq n_0$ for any fixed $n_0$.
If we define $d_n=E_n-E_{n-1}$ for $n\ge 0$
(with the convention that $E_{-1}=0$), for $j\geq1$ we have
$$
\Delta_j=E_{2^j}-E_{2^{j-1}}=\sum_{k=0}^{2^j}d_k-\sum_{k=0}^{2^{j-1}}d_k=
\sum_{k=2^{j-1}+1}^{2^j}d_k.
$$
Consider, for each $n\geq 5$, $J_n$ the unique integer such that
$2^{J_n}<n\leq 2^{J_n+1}$. Then
$$
\sigma_n-\sigma_{n-1}=\frac{1}{n(n+1)}\sum_{j=0}^n jd_j
=\frac{1}{n(n+1)}\Big(d_1+2d_2+\sum_{k=2}^{J_n}\sum_{j=2^{k-1}+1}^{2^k}
jd_j +\sum_{k=2^{J_n}+1}^{n} jd_j\Big).
$$
Now, for each $k$, $2\leq k\leq J_n$,
\begin{eqnarray*}
\sum_{j=2^{k-1}+1}^{2^k} jd_j
&=&
2^k\Delta_k-\sum_{j=2^{k-1}+1}^{2^k} (2^k-j)d_j\\
&=&
2^k\Delta_k-\sum_{j=2^{k-1}+1}^{2^k-1} \sum_{i=j}^{2^k-1}d_j
=
2^k\Delta_k-\sum_{i=2^{k-1}+1}^{2^k-1} E_i(\Delta_k).
\end{eqnarray*}
We can treat the rest of the terms in a similar way, and then we get
\begin{eqnarray*}
\lefteqn{\sigma_n-\sigma_{n-1}=
\frac{1}{n(n+1)}\Big[d_1+2d_2+\sum_{k=2}^{J_n}
\Big(2^k\Delta_k-\sum_{i=2^{k-1}+1}^{2^k-1} E_i(\Delta_k)\Big)+}\\
& &\hspace{3.5cm}
+nE_n(\Delta_{J_n+1})
-\sum_{k=2^{J_n}+1}^{n-1}E_k(\Delta_{J_n+1})\Big].
\end{eqnarray*}
Thus
\begin{eqnarray*}
\lefteqn{\Big(\sum_{n=5}^\infty n^{q-1}\|(\sigma_n-\sigma_{n-1})
f\|_\B^q\Big)^{1/q} \leq \|d_1f+2d_2f\|_\B\Big(\sum_{n=5}^\infty
\frac{n^{q-1}}{n^q(n+1)^q}\Big)^{1/q}}\\
&+ &
\Big(\sum_{n=5}^\infty
\frac{n^{q-1}}{n^q(n+1)^q}
\Big\|\sum_{k=2}^{J_n}2^k\Delta_kf\Big\|_\B^q\Big)^{1/q}
+ \Big(\sum_{n=5}^\infty \frac{n^{q-1}}{n^q(n+1)^q}
\Big\|\sum_{k=2}^{J_n}\sum_{i=2^{k-1}+1}^{2^k-1} E_i(\Delta_kf)
\Big\|_\B^q\Big)^{1/q}\\
&+& \Big(\sum_{n=5}^\infty \frac{n^{q-1}}{n^q(n+1)^q}
\|nE_n(\Delta_{J_n+1}f)\|_\B^q\Big)^{1/q}\\
&+& \Big(\sum_{n=5}^\infty\frac{n^{q-1}}{n^q(n+1)^q}
\Big\|\sum_{k=2^{J_n}+1}^{n-1}E_k(\Delta_{J_n+1}f)\Big\|_\B^q\Big)^{1/q}\\
&=& C\|d_1f+2d_2f\|_\B+I+II+III+IV.
\end{eqnarray*}
Using $ \big|\sum_{i=1}^n2^ia_i\big|^q\leq
2^{n(q-1)}\sum_{i=1}^n2^i|a_i|^q$,
 we have that
\begin{eqnarray*}
I^q
&\leq&
\sum_{n=5}^\infty \frac{1}{n^{q+1}}2^{J_n(q-1)}
\sum_{k=2}^{J_n}2^k\|\Delta_kf\|_\B^q
\leq
\sum_{n=5}^\infty \frac{1}{n^2}
\sum_{k=2}^{J_n}2^k\|\Delta_kf\|_\B^q\\
&\leq&
\sum_{k=2}^\infty\|\Delta_kf\|_\B^q.
\end{eqnarray*}
Since  $\B$ is of martingale cotype $q$, $\|I\|_{L^p}\leq
C_{p,q,\B}\;\|f\|_{L^p_\B}$.

In order to handle the second term, let
us call $k(j)=k$ when $2^{k-1}\leq j<2^k$. Then,
\begin{eqnarray*}
II^q &\leq& \sum_{n=5}^\infty \frac{1}{n^{q+1}}
\Big\|\sum_{j=3}^{2^{J_n}-1} E_j(\Delta_{k(j)}f)\Big\|_\B^q \leq
\sum_{n=5}^\infty \frac{1}{n^{q+1}}2^{J_n(q-1)}
\sum_{j=3}^{2^{J_n}-1} \|E_j(\Delta_{k(j)}f)\|_\B^q\\
&\leq& \sum_{n=5}^\infty \frac{1}{n^2} \sum_{j=3}^{n-1}
\|E_j(\Delta_{k(j)}f)\|_\B^q \leq \sum_{j=3}^\infty\frac{1}{j+1}
\|E_j(\Delta_{k(j)}f)\|_\B^q.
\end{eqnarray*}
Using (\ref{p inequality}) and the martingale cotype $q$ of $\B$,
we obtain that
\begin{eqnarray*}
\|II\|_{L^p} &\leq& \Big\|\Big(\sum_{j=3}^\infty\frac{1}{j+1}
\|E_j(\Delta_{k(j)}f)
\|_\B^q\Big)^{1/q}\Big\|_{L^p}\\
&\leq& C_{p,q}\;\Big\|\Big(\sum_{k=2}^\infty \|\Delta_{k}f\|_\B^q
\sum_{j=2^{k-1}+1}^{2^k}\frac{1}{j+1}\Big)^{1/q}\Big\|_{L^p} \leq
C_{p,q,\B}\;\|f\|_{L^p_\B}.
\end{eqnarray*}
 Analogously, one can show that
$\|III\|_{L^p}+\|IV\|_{L^p}\leq C_{p,q}\|f\|_{L^p_\B}$. We leave
the details to the reader. Thus the lemma is proved.\cqd

\medskip

We shall also need the following result due to Rota, see \cite
[Chapter IV] {st:}. Let $Q$ be a linear operator on
$L^p(\Omega,{\mathcal A} ,d\mu)$ satisfying
the conditions\\
\noindent{\sf i)} $\|Q\|_{L^p\to L^p}\leq1$ for every $p\in[1,\infty]$,\\
\noindent{\sf ii)} $Q=Q^*$ in $L^2$,\\
\noindent{\sf iii)} $Qf\geq 0$ for every $f\geq 0$,\\
\noindent{\sf iv)} $Q1=1$.

\begin{theorem}\label{rota}
For any  $Q$ as above, there exist a measure space $(M,\F,dm)$, a
decreasing collection of $\sigma$-algebras
$\dots\subset\F_{n+1}\subset\F_n\subset\dots
\subset\F_1\subset\F_0\subset\F$, and another $\sigma$-algebra
$\hat\F\subset\F$
such that\\
\noindent {\sf a)} there exists an isomorphism $i:(\Omega,{\mathcal
A},d\mu)\rightarrow(M,\hat\F,dm)$ (which induces an isomorphism
between $L^p$ spaces, also denoted by $i$,
$i(f)(m)=f(i^{-1}m)$),\\
\noindent {\sf b)} for every $f\in L^p(M,\hat\F,dm)$, we have
$$Q^{2n}(i^{-1}f)(x)=\hat E(E_n(f))(i(x)),\quad x\in \Omega$$
where $\hat E(f)=E(f|\hat\F)$ and $E_n(f)=E(f|\F_n)$.
\end{theorem}
This theorem holds in the scalar valued case. For the vector
valued case, the validity of the second statement is a consequence
of the fact that all operators in consideration extend to
contractions on $\B$-valued $L^p$-spaces. Indeed, the linearity
implies that the formula in the statement {\sf b)} above holds for
all $\B$-valued simple functions, and so for all $\B$-valued
$p$-integrable functions.

\medskip

\noindent{\sc Proof of Theorem \ref{theorem on the means}}.
Observe that it is enough to prove
$$
\Big\|\Big(\int_a^b \Big\|t\frac{\partial M_tf}{\partial
t}\Big\|^q_\B \,\frac{dt}{t}\Big)^{1/q}\Big\|_{L^p(\Omega)} \leq
C_{p,q,\B} \|f\|_{L^p_\B(\Omega)}
$$
for any $0<a<b<\infty$, and also that it is enough if we restrict
ourselves to functions $f$ in the algebraic tensor product
$\B\otimes L^p(\Omega)$. Take then $f=\sum_{k=1}^Kv_k\varphi_k$.
By the results in \cite{st:}, see the lemma in p.72 and its proof,
it is not difficult to observe that for every $t_0\in(0,\infty)$,
there exists $\varepsilon_0>0$ such that
\begin{equation}
\label{power series} \Tt_tf(x)=\sum_{j=0}^\infty f_j(x)(t-t_0)^j
\end{equation}
for $t\in (t_0-\varepsilon_0,t_0+\varepsilon_0)$ and almost every
$x$ , where $\sum\|f_j\|_{L^p_\B}\;\varepsilon_0^j<\infty$ and
where $f_j$ depend on $t_0$. Since we can cover $(a,b)$ with a
finite collection of such intervals, we can split $(a,b)$ into a
finite collection of subintervals $(a_i,b_i)$ of $(a,b)$ in which
a expression like \eqref{power series} holds for a fixed $t_0$
(and therefore, with the same $f_j$) for every $t\in (a_i,b_i)$.
Then, splitting the integral between $a$ and $b$ into the
integrals corresponding to such subintervals, we can handle all
the functions appearing in the integral as power series with
vector valued coefficients. In these circumstances, we can replace
the integral by Riemann sums, and all derivatives
by difference quotients. 
The first step is choosing $\varepsilon$ small. Then, we
approximate the integral as follows:
$$
\Big\|\Big(\int_a^b \Big\|t\frac{\partial M_tf}{\partial
t}\Big\|^q_\B \,\frac{dt}{t}\Big)^{1/q}\Big\|_{L^p} \sim
\Big\|\Big(\sum_{n=n_0}^{n_1}(n\varepsilon)^{q-1}
\Big\|\frac{\partial M_tf}{\partial
t}\Big|_{t=n\varepsilon}\Big\|^q_\B\varepsilon
\Big)^{1/q}\Big\|_{L^p},
$$
where the sign $\sim$ means that the difference term goes to zero
as $\varepsilon\to 0$. The next step is substituting the partial
derivative inside the sum by the difference quotient
$$
\frac{M_{(n+1)\varepsilon}f-M_{n\varepsilon}f}{\varepsilon}
=
\frac1\varepsilon\frac{1}{(n+1)\varepsilon}\int_0^{(n+1)\varepsilon}\Tt_sf\,ds-
\frac1\varepsilon\frac{1}{n\varepsilon}\int_0^{n\varepsilon}\Tt_sf\,ds,
$$
and then each of the integrals by its Riemann sums,
getting then that
\begin{eqnarray*}
\Big\|\Big(\int_a^b \Big\|t\frac{\partial M_tf}{\partial
t}\Big\|^q_\B \,\frac{dt}{t}\Big)^{1/q}\Big\|_{L^p} &\sim&
\Big\|\Big(\sum_{n=n_0}^{n_1}n^{q-1}\Big\|\frac{1}{n+1}
\sum_{j=0}^n\Tt_{j\varepsilon}f-
\frac{1}{n}\sum_{j=0}^{n-1}\Tt_{j\varepsilon}f
\Big\|^q_\B\Big)^{1/q}\Big\|_{L^p}\\
&=&
\Big\|\Big(\sum_{n=n_0}^{n_1}n^{q-1}\|\tilde\sigma_nf-\tilde\sigma_{n-1}f
\|^q_\B\Big)^{1/q}\Big\|_{L^p},
\end{eqnarray*}
where
$\tilde\sigma_nf=\frac{1}{n+1}\sum_{j=0}^n\Tt_{j\varepsilon}f$.
Now, observe that by our hypothesis, $\Tt_{\varepsilon/2}$
satisfies assumptions {\sf i)}-{\sf iv)} of Rota's theorem and
$\Tt_{n\varepsilon}f=\hat E(E_n(f))$. Hence, $\tilde\sigma_nf=\hat
E(\sigma_nf)$ where $\sigma_n$ is as in Lemma \ref{cesaro}.
Therefore, by the properties of conditional expectation and
Lemma \ref{cesaro}, we get
\begin{eqnarray*}
 \Big\|\Big(\int_a^b \Big\|t\frac{\partial M_tf}
 {\partial t}\Big\|^q_\B \,\frac{dt}{t}\Big)^{1/q}\Big\|_{L^p}
 &\sim&\Big\|\Big(\sum_{n=n_0}^{n_1} n^{q-1}\|\hat E(\sigma_nf)-\hat
 E(\sigma_{n-1}f)\|^q_\B \Big)^{1/q}\Big\|_{L^p}\\
 &\leq&\Big\|\big\{n^{1-1/q}(\sigma_nf-\sigma_{n-1}f)\big\}\Big\|
 _{L^p_{\ell^q_\B}}\!\!\\
 &=& \Big\|\Big(\sum_{n=1}^\infty n^{q-1}\|\sigma_nf-
 \sigma_{n-1}f\|^q_\B\Big)^{1/q}\Big\|_{L^p}\\
 &\leq&\!\!C_{p,q,\B}\; \|f\|_{L^p_\B}\ .
\end{eqnarray*}
Therefore, we have achieved the proof of Theorem \ref{theorem on
the means}.\cqd

\medskip

The following lemma says that the boundedness of
 $G_qf=\|Tf\|_{L_\B^q((0,1),\;dr/(1-r))}$
from $L_\B^p(\T)$ in $L^p(\T)$ is equivalent to the boundedness of
the operator $T$ when the kernel is restricted to values of $r$
close to one and $\theta$ close to zero.

\begin{lemma}\label{restriction g in the torus 1}
Let $\B$ be a Banach space and $p,q\in(1,\infty)$. Let $\delta>0$
(close to 0). Then there is a constant $C_\delta$ (depending only
on $\delta$) such that for any $f\in L^p_\B(\T)$
$$
\Big\|\big[(1-r)\frac{\partial P_r}{\partial
r}\ind_{(0,1-\delta)}(r)\big]
*f\Big\|_{L^p_{L^q_\B((0,1),\frac{dr}{1-r})}(\T)} \leq C_\delta
\|f\|_{L^p_\B(\T)}\ ,
$$
and
$$
\Big\|\big[(1-r)\frac{\partial P_r}{\partial
r}\ind_{(1-\delta,1)}(r) \ind_{(-\delta, \delta)}(\theta)\big]
*f\Big\|_{L^p_{L^q_\B((0,1),\frac{dr}{1-r})}(\T)} \leq C_\delta
\|f\|_{L^p_\B(\T)}.
$$
\end{lemma}

\noindent{\sc Proof}. The proof is very easy. We  show only the
first inequality. Its left hand side is a convolution of the
$\B$-valued function $f$ with an
$L^q((0,1),\frac{dr}{1-r})$-valued function. Therefore it is
enough to prove that the latter is in
$L^1_{L^q((0,1),\frac{dr}{1-r})}$, namely, we have to prove
$$
\int_0^{2\pi} \Big\|(1-r)\frac{\partial P_r(\gamma)}{\partial
r}\ind_{(0,1-\delta)}(r) \Big\|_{L^q((0,1),\frac{dr}{1-r})}
\,d\gamma <\infty.
$$
But this follows immediately if we observe that
\begin{eqnarray*}
\lefteqn{\Big\|(1-r)\frac{\partial P_r(\gamma)}{\partial
r}\ind_{(0,1-\delta)}(r)
\Big\|_{L^q((0,1),\frac{dr}{1-r})}^q}\\
&=&
\int_0^{1-\delta}\Big|(1-r)\frac{2(1-r)^2-2(r^2+1)\sin^2(\gamma/2)}
{((1-r)^2+2r\sin^2(\gamma/2))^2}\Big|^q\, \frac{dr}{1-r}\leq
C_\delta^q.
\end{eqnarray*}
Hence the lemma is proved\cqd

\medskip
The following easy lemma is proved in a similar way as in \cite
[p.49]{st:}.

\begin{lemma}\label{poisson-mean}
Let $\B$ be a Banach space and $p, q\in(1,\infty)$. Then for any
$f\in L^p(\Omega)\otimes \B$ we have
 $$\Big(\int_0^\infty\Big\|t\frac{\partial
 \Pt_tf}{\partial t}\Big\|_\B^q\,\frac{dt}{t}\Big)^{1/q} \leq
 C_0\Big(\int_0^\infty \Big\|t\frac{\partial M_tf}{\partial
 t}\Big\|_\B^q\,\frac{dt}{t}\Big)^{1/q}\ ,$$
where $C_0$ is an absolute constant.
\end{lemma}

\noindent{\sc Proof}.  Call
$\varphi(s)=\frac{1}{2\sqrt{\pi}}\frac{e^{-1/4s}}{s^{3/2}}$. Using
integration by parts we have
$$
 \Pt_t
 = \frac{1}{t^2}\int_0^\infty\varphi (\frac{s}{t^2})
 \big(\frac{\partial}{\partial s}sM_s\big)\,ds
 =-\int_0^\infty \frac{s}{t^4}
 \varphi'\big(\frac{s}{t^2}\big)\,M_s\,ds
 =-\int_0^\infty s\,
 \varphi'\big(s)\,M_{t^{2}s}\,ds.
$$
Therefore (with $M'_s=\frac{\pa M_s}{\pa s}$),
\begin{equation}\label{poisson-mean1}
 t\frac{\partial}{\partial t}\Pt_t=
 -2\int_0^\infty t^{2}\,s^2\varphi'(s)  M'_{t^{2}s}\,ds
 =-2\int_0^\infty s\,\varphi'(s)\big[t^{2}s\, M'_{t^{2}s}\big]ds\ .
\end{equation}
Thus
\begin{eqnarray*}
 \Big[\int_0^\infty\Big\|t\frac{\partial \Pt_tf}{\partial
 t}\Big\|_\B^q\,\frac{dt}{t}\Big]^{1/q}
 &\leq& 2\int_0^\infty s\,|\varphi'(s)|
 \Big[\int_0^\infty\big\|t^{2}s\, M'_{t^{2}s}f\big\|_\B^q\,
 \frac{dt}{t}\Big]^{1/q}ds\\
 &=&2^{1-1/q}K \Big[\int_0^\infty\big\|t M'_tf\big\|_\B^q\,
 \frac{dt}{t}\Big]^{1/q}\;,
\end{eqnarray*}
where
$$K=\int_0^\infty s\,|\varphi'(s)|\,ds.$$
Hence the lemma is proved. \cqd

\medskip

Now we are well prepared for the proofs of Theorems \ref{lusin
semigroups} and \ref{lusin type semigroups}.

\medskip

\noindent{\sc Proof of Theorem \ref{lusin semigroups}}. {\sf i)}
$\Rightarrow$ {\sf ii)}. This is an immediate consequence of
Theorem \ref{theorem on the means} and Lemma \ref{poisson-mean}.

{\sf ii)}$\Rightarrow$ {\sf i)}. We shall prove that the operator
$f\mapsto G^1_q(f)$ is bounded from $L_\B^p(\T)$ to $L^p(\T)$ for
$p\in(1,\infty)$. Recall that
$$
 G^1_q(f)(z)=\Big(\int_0^1 (1-r)^q\|\frac{\partial P_r}{\partial
 r}*f(z)\|_{\B}^q\,\frac{dr}{1-r}\Big)^{1/q}\ ,\quad z\in\T.
$$
By \cite{xu}, this is equivalent to the
martingale cotype $q$ of $\B$.  Observe that if in the Poisson
kernel $P_r$, $0<r<1$, we change the parameter according to
$r=e^{-t}$, we obtain the kernel $\widetilde P_t$ of the Poisson
semigroup subordinated to the heat semigroup in the torus. Fix a
$\delta\in(0,1)$ (very close to 1). By the same change of
parameter and the fact that for any $t\in(0,-\log \delta)$,
$\frac{1-e^{-t}}{e^{-t}}\sim t$, then we have
\begin{eqnarray*}
\lefteqn{\int_{\delta}^1\big\|(1-r)\frac{\partial P_r}{\partial r}
*f(\theta)\big\|^q\,\frac{dr}{1-r} =
\int_0^{-\log\delta}\big\|\frac{1-e^{-t}}{e^{-t}}\frac{\partial
\widetilde P_t}{\partial t}*f(\theta)\big\|^q\,
\frac{e^{-t}dt}{1-e^{-t}}}\\
&\hspace{1cm}&\leq
C_{\delta,q}^q\int_0^{-\log\delta}\big\|t\frac{\partial \widetilde
P_t}{\partial t}*f(\theta)\big\|^q\, \frac{dt}{t} \leq
C_{\delta,q}^q\int_0^\infty\big\|t\frac{\partial \widetilde
P_t}{\partial t}*f(\theta)\big\|^q\, \frac{dt}{t}.
\end{eqnarray*}
Therefore, by hypothesis {\sf ii)}, we have that
$$
\Big\|\big[(1-r)\frac{\partial P_r}{\partial
r}\ind_{(\delta,1)}(r)\big]
*f\Big\|_{L^p_{L^q_\B((0,1),\frac{dr}{1-r})}(\T)} \leq
C'_{\delta,q} \|f\|_{L^p_\B(\T)}.
$$
Then by Lemma \ref{restriction g in the torus 1} we get
$$
\big\|G^1_q(f)\big\|_{L^p(\T)} \leq C^{''} \|f\|_{L^p_\B(\T)}.
$$
By  \cite{xu}, this implies that $\B$ is of Lusin cotype $q$, and
so of martingale cotype $q$ too. Thus the proof of Theorem
\ref{lusin semigroups} is finished. \cqd

\medskip

\noindent{\sc Proof of Theorem \ref{lusin type semigroups}}. {\sf
i)} $\Rightarrow$ {\sf ii)}. Write the spectral decomposition of
the semigroup $\{\Pt_t\}_{t\geq 0}$: for any $f\in L^2(\Omega)$
$$
\Pt_tf=\int_0^\infty e^{-\lambda t}de_\lambda f,
$$
where $\{e_\lambda\}$ is a resolution of the identity. Thus
$$
\frac{\partial \Pt_tf}{\partial t}=-\int_{0+}^\infty \lambda
e^{-\lambda t}de_\lambda f.
$$
It is easy to deduce from this formula that for any $f, g\in
L^2(\Omega)$ (recalling that $F$ is the projection on the fix
point subspace of $\{\Pt_t\}_{t\geq 0}$)
\begin{equation}
\label{parseval}
 \int_\Omega (f-F(f))(g-F(g))d\mu=4
 \int_\Omega \int_0^\infty\big[t\frac{\partial \Pt_tf}{\partial
 t}\big]\big[t\frac{\partial \Pt_tg}{\partial
 t}\big]\,\frac{dt}{t}\,d\mu.
\end{equation}

Now we use duality. Fix two functions $f\in L_\B^p(\Omega)$ and
$g\in L_{\B^*}^{p'}(\Omega)$, where $p'$ denotes the conjugate
index of $p$ . Without loss of generality we may assume that $f$
and $g$ are in the algebraic tensor products $\big(L^p(\Omega)\cap
L^2(\Omega)\big)\otimes\B$ and $\big(L^{p'}(\Omega)\cap
L^2(\Omega)\big)\otimes\B^*$, respectively. With $\langle\,
,\,\rangle$ denoting the duality between $\B$ and $\B^*$, we have
$$
 \int_\Omega \langle f,\; g\rangle d\mu
 =\int_\Omega \langle F(f),\; F(g)\rangle d\mu
 +\int_\Omega \langle f-F(f),\; g-F(g)\rangle\, d\mu.
$$
The first term on the right is easy to be estimated:
$$
 \Big|\int_\Omega \langle F(f),\; F(g)\rangle d\mu\Big|
 \le \|F(f)\|_{L_\B^p(\Omega)}\;\|F(g)\|_{L_{\B^*}^{p'}(\Omega)}
 \le \|F(f)\|_{L_\B^p(\Omega)}\;\|g\|_{L_{\B^*}^{p'}(\Omega)}.
$$
For the second one, by  (\ref{parseval}) and H\"older's inequality
\begin{eqnarray*}
 \Big|\int_\Omega\langle f-F(f),\; g-F(g)\rangle d\mu\Big|
 &=&4\Big|\int_\Omega\int_0^\infty\langle t\frac{\partial\Pt_tf}{\partial
 t},\;t\frac{\partial\Pt_tg}{\partial
 t}\rangle\,\frac{dt}{t}\,d\mu\Big|\\
 &\le& 4\int_\Omega\int_0^\infty\|t\frac{\partial\Pt_tf}{\partial
 t}\|\;\|t\frac{\partial \Pt_tg}{\partial
 t}\|\,\frac{dt}{t}\,d\mu\\
 &\le&4\|{\mathfrak G}_q(f)\|_{L^p(\Omega)}
 \|{\mathfrak G}_{q'}(g)\|_{L^{p'}(\Omega)}\ .
\end{eqnarray*}
Now since $\B$ is of martingale type $q$, $\B^*$ is of martingale
cotype $q'$. Thus by Theorem \ref{lusin semigroups},
$$
 \|{\mathfrak G}_{q'}(g)\|_{L^{p'}(\Omega)}
 \le C\;\|g\|_{L_{\B^*}^{p'}(\Omega)}\ .
$$
Combining the preceding inequalities, we get
$$
 \Big|\int_\Omega \langle f,\ g\rangle d\mu\Big|
 \le \big(\|F(f)\|_{L_\B^p(\Omega)} +
 C\|{\mathfrak G}_q(f)\|_{L^p(\Omega)}\big)
 \|g\|_{L_{\B^*}^{p'}(\Omega)}\ .
$$
which gives {\sf ii)}, taking the supremum over all $g$ as above such that
$\|g\|_{L_{\B^*}^{p'}(\Omega)}\le 1$.

{\sf ii)} $\Rightarrow$ {\sf i)}. As in the corresponding proof of
Theorem \ref{lusin semigroups}, we use again the Poisson semigroup
on the torus. We keep the notations introduced there. Recall that
$\widetilde P_t=P_{e^{-t}}$. By the calculations done there,
$$
 \int_{\delta}^1\big\|(1-r)\frac{\partial P_r}{\partial r}
 *f(\theta)\big\|^q\,\frac{dr}{1-r} \approx
 \int_0^{-\log\delta}\big\|t\frac{\partial \widetilde P_t}{\partial
 t}*f(\theta)\big\|^q\, \frac{dt}{t},
$$
where the equivalence constants depend only on $\delta$ and $q$.
On the other hand, on the interval $(0, \delta)$, we have
\begin{eqnarray*}
 \int_0^{\delta}\big\|(1-r)\frac{\partial P_r}{\partial r}
 *f(\theta)\big\|^q\,\frac{dr}{1-r}
 &=&\int_{-\log\delta}^\infty\big\|\frac{1-e^{-t}}{e^{-t}}
 \frac{\partial\widetilde P_t}{\partial t}*f(\theta)\big\|^q\,
 \frac{e^{-t}dt}{1-e^{-t}}\\
 &\ge& C_{\delta,q}^q\int_{-\log\delta}^\infty e^{(q-1)t}
 \big\|\frac{\partial \widetilde P_t}
 {\partial t}*f(\theta)\big\|^q\,dt\\
 &\ge&(C'_{\delta,q})^q\int_{-\log\delta}^\infty
 t^{q-1}\big\|\frac{\partial \widetilde P_t}{\partial
 t}*f(\theta)\big\|^q\, dt .
\end{eqnarray*}
Therefore,
$$
 \int_{0}^\infty t^{q-1}\big\|
 \frac{\partial\widetilde P_t}{\partial t}*f(\theta)\big\|^q\,dt
 \le C_{\delta,q}^q \int_{0}^1\big\|(1-r)\frac{\partial P_r}{\partial r}
 *f(\theta)\big\|^q\,\frac{dr}{1-r}\ .
$$
Thus by hypothesis {\sf ii)},
$$
 \|f\|_{L_\B^p(\T)}\le C_{\delta, q}\big(\|\hat f(0)\|_\B+
 \|G^1_q(f)\|_{L^p(\T)}\big)\ .
$$
Hence by \cite{xu}, $\B$ is of Lusin type $q$, and so of
martingale type $q$ too. \cqd

\medskip

We end this section with some remarks and questions.

\begin{remark} Checking back the proofs above of Theorems
\ref{lusin semigroups} and \ref{lusin type semigroups}, we see
that under the condition that $\B$ is of martingale cotype $q$
(resp. martingale type $q$), (\ref{cotype semigroup}) (resp.
(\ref{type semigroup})$\,$) is true for more  general semigroups
$\{\Pt_t\}$ associated to a symmetric diffusion semigroup other
than the subordinated Poisson semigroup given by
(\ref{subordinated semigroup}). ($\{\Pt_t\}$ needs not be even a
semigroup for the validity of  (\ref{cotype semigroup}).) What we
need is that $\{\Pt_t\}$ is defined by
$$
 \Pt_t = \int_0^\infty
 \varphi(s)\Tt_{t^\beta s}\,ds,$$
where $\beta$ is a non zero real number, and where $\varphi$ is a
derivable function on $\R_+$ such that both $\varphi$ and
$t\,\varphi'$ are integrable on $\R_+$  and such that the two
limits $\lim_{t\to0}t\varphi(t)$ and
$\lim_{t\to\infty}t\varphi(t)$ exist. For instance, this is the
case when the infinitesimal generator of $\{\Pt_t\}$ is
$-(-A)^\alpha$ with $0<\alpha<1$, where $A$ is the infinitesimal
generator of a symmetric diffusion semigroup $\{\Tt_t\}$. Indeed,
by \cite [IX] {yosido}, $\{\Pt_t\}$ can be represented as above
with $\beta=1/\alpha$ and $\varphi$ given by
$$
 \varphi(s)=\int_0^\infty \exp\big[st\cos\theta -
 t^\alpha\cos(\alpha\theta)\big]\times\sin\big[st\sin\theta -
 t^\alpha\cos(\alpha\theta)+\theta\big]dt,
$$
where $\theta$ can be any number in $[\pi/2,\ \pi]$.
\end{remark}

\begin{remark} The proof of {\sf i)} $\Rightarrow$ {\sf ii)}
in Theorem \ref{lusin type semigroups} implicitly shows the
following: Given a Banach space $\B$ and $p, q\in(1,\infty)$, if
(\ref{cotype semigroup}) holds for a given subordinated semigroup
$\{\Pt_t\}$, then (\ref{type semigroup}) holds for the same
semigroup with $\B$ and $p, q\in(1,\infty)$ replaced by $\B^*$ and
$p', q'\in(1,\infty)$, respectively. The converse to this latter
statement is also true. This will be the objective of the next
section. With this converse, we can prove {\sf ii)} $\Rightarrow$
{\sf i)} in Theorem \ref{lusin type semigroups} directly by
duality and Theorem \ref{lusin semigroups} without using the
Poisson semigroup on $\T$. Note that such an approach is
inevitable when one wishes to study the duality between
(\ref{cotype semigroup}) and (\ref{type semigroup}) for an
individual semigroup.
\end{remark}

All semigroups considered in this paper are Markovian, that is,
$\Tt_t1=1$. We do not know whether Theorems  \ref{lusin
semigroups} and  \ref{lusin type semigroups} still hold for
symmetric sub-Markovian semigroups (which are those satisfying
(\ref {semigroup 1}) - (\ref {semigroup 3}) except the
Markovianity).

\medskip

\noindent{\bf Problem 1.} {\em Let $\B$ be a Banach space of
martingale cotype $q$. For which semigroups $\{\Pt_t\}$  is the
corresponding $g$-function mapping $f\mapsto {\mathfrak G}_q(f)$
of weak type (1,1)?}

\medskip

We shall see later that the answer is positive for the usual
Poisson semigroup and the subordinated Poisson Ornstein-Uhlenbeck
semigroup on $\R^n$. We shall also show that (\ref{cotype
semigroup}) (resp. (\ref{type semigroup})) for one of these
Poisson semigroups implies the martingale cotype $q$ (resp.
martingale type $q$) of $\B$, like for the Poisson semigroup on
the torus.

In general, it would be interesting to find conditions on a given
semigroup $\{\Pt_t\}$ which guarantee that the validity of
(\ref{cotype semigroup}) (resp. (\ref{type semigroup})) for
$\{\Pt_t\}$ implies martingale cotype $q$ (resp. martingale cotype
$q$ ).

\medskip

We state another problem about (\ref{cotype semigroup}) for any
symmetric diffusion semigroup (not necessarily subordinated to
another one).

\medskip

\noindent{\bf Problem 2.} {\em Let $\B$ be a Banach space of
martingale cotype $q$ (resp. martingale type $q$) and
$p\in(1,\infty)$. Does (\ref{cotype semigroup}) (resp. (\ref{type
semigroup})) hold for any symmetric diffusion semigroup
$\big\{\Pt_t\big\}_{t\ge0}$?}

\medskip

Problem 2 has an affirmative solution when $\B$ is further a
Banach lattice. Let us consider only the cotype case. It is well
known that a Banach lattice $\B$ is of martingale cotype $q$ (with
$2<q<\infty$) iff $\B$ is $q$-concave and $p$-convex for some
$p>1$. For $q=2$, the \lq\lq if \rq\rq\ part is still true; the
\lq\lq only if\rq\rq\ part admits only a weaker form: $\B$ is
$r$-concave and $p$-convex for some $p>1$ and for any $r>q=2$. See
\cite {LT}.  Let $\B$ be a $q$-concave and $p$-convex Banach
lattice with $1<p\le 2\le q<\infty$. Then by \cite {pi:79}, $\B$
can be written as a complex interpolation space between a Hilbert
space and another lattice, i.e., there are a Hilbert space $H$ and
a Banach lattice $\B_0$ such that
 $$\B=\big(H,\ \B_0\big)_{\theta}\ .$$
Moreover, $\B_0$ is $q_0$-concave and $p_0$-convex with $p_0>1$
and $q_0$ satisfying $1/q=\theta/2 + (1-\theta)/q_0$. Now given a
symmetric diffusion semigroup $\{\Tt_t\}$, following \cite
[p.116-119] {st:}, we consider the fractional averages of
$\{\Tt_t\}$:
 $$M_t^\alpha(f)=\frac{t^{-\alpha}}{\Gamma(\alpha)}
 \int_0^t(t-s)^{\alpha -1}\;T_sf\,ds.$$
$M_t^\alpha$ is well defined for $\alpha\in \Comp$ with ${\rm
Re}(\alpha)>0$, and is continued analytically into the whole
complex plane. Note that $M_t^1$ is the usual average $M_t$ and
$M_t^0=T_t$. By \cite {st:}, for all $\alpha\in\Comp$ and $f\in
L^p_H(\Omega)$ ($1<p<\infty$), we have
 $$\Big\|\big[\int_0^\infty\big\|t\,\frac{\pa M_t^\alpha f}{\pa t}
 \big\|_H^2\,\frac{dt}{t}\big]^{1/2}\Big\|_{L^p(\Omega)}\le
 C_{p,\alpha}\;\|f\|_{L^p_H(\Omega)}\ .$$
(This is proved in \cite{st:} for the scalar valued case; but the
the same arguments work as well for Hilbert space valued
functions.) On the other hand, using Theorem \ref {theorem on the
means}, one can easily check that for any $\alpha\in\Comp$ with
${\rm Re}(\alpha)>1$ and $f\in L^p_{\B_0}(\Omega)$ ($1<p<\infty$)
 $$\Big\|\big[\int_0^\infty\big\|t\,\frac{\pa M_t^\alpha f}{\pa t}
 \big\|_{\B_0}^{q_0}\,\frac{dt}{t}\big]^{1/q_0}\Big\|_{L^p(\Omega)}\le
 C_{p,\B_0,\alpha}\;\|f\|_{L^p_{\B_0}(\Omega)}\ .$$
Then interpolating these inequalities, we deduce that for any
$f\in L^p_{\B}(\Omega)$ ($1<p<\infty$)
 $$\Big\|\big[\int_0^\infty\big\|t\,\frac{\pa M_t^0 f}{\pa t}
 \big\|_\B^{q}\,\frac{dt}{t}\big]^{1/q}\Big\|_{L^p(\Omega)}\le
 C_{p,\B}\;\|f\|_{L^p_{\B}(\Omega)}\ .$$
This is the desired inequality on $T_t$ (recalling that
$T_t=M_t^0$).


\section{Duality}

Throughout this section  $\{\Tt_t\}_{t\geq 0}$ will be a fixed
symmetric diffusion semigroup defined on $L^p(\Omega,d\mu)$, and
$\{\Pt_t\}_{t\geq 0}$ its subordinated Poisson semigroup. We shall
keep all notations introduced in the previous section for these
semigroups. In particular, $F$ is the contractive projection from
$L^p(\Omega)$ (also from $L_\B^p(\Omega)$) onto the fix point
subspace of $\{\Pt_t\}_{t\geq 0}$. The following is the main
result of this section.

\begin{theorem}\label{type-cotype semigroup}
Let $\B$ be a Banach space and $1<p, q<\infty$. Then the following
statements are equivalent:
\begin{itemize}
\item[{\sf i)}] There is a constant $C>0$ such that
$$
\|{\mathfrak G}_qf\|_{L^p(\Omega)}\leq C\|f\|_{L^p_\B(\Omega)},
\quad \forall\;f\in L^p_\B(\Omega).
$$
\item[{\sf ii)}] There is a constant $C>0$ such that
$$
\|g\|_{L^{p'}_{\B^*}(\Omega)}\le
C\big(\|F(g)\|_{L^{p'}_{\B^*}(\Omega)} + \|{\mathfrak
G}_{q'}g\|_{L^{p'}(\Omega)}\big), \quad \forall\;g\in
L^{p'}_{\B^*}(\Omega).
$$
\end{itemize}
\end{theorem}

The proof of the implication {\sf i)} $\Rightarrow$ {\sf ii)} of
Theorem \ref{lusin type semigroups} shows in fact {\sf i)}
$\Rightarrow$ {\sf ii)} in the theorem above. The inverse
implication needs much more effort (as usually in such a
situation). Let $\AB=L_\B^q(\R_+, \frac{dt}{t})$. An element $h$
in $L^p_{\AB}(\Omega)$ is a function of two variables $x\in\Omega$
and $t\in\R_+$, i.e., $h: (x, t)\mapsto h_t(x)$. The key to the
implication {\sf ii)} $\Rightarrow$ {\sf i)} above is the
existence of a bounded projection from $L^p_{\AB}(\Omega)$ onto
the subspace of all functions $h$ which can be written as
$\displaystyle h_t(x)=t\frac{\pa \Pt_tf}{\pa t}(x)$ for some
function $f$ on $\Omega$. Formally, the desired projection is
given by $\displaystyle h\mapsto t\frac{\pa \Pt_t(Qh)}{\pa t}(x)$,
where $Qh$ is defined by
 \begin{equation}\label{Q}
 Qh(x)=\int_0^\infty t\frac{\pa \Pt_th_t}{\pa t}(x)\; \frac{dt}{t}\
 ,\quad x\in\Omega.
 \end{equation}
Note that $Q(h)$ is well-defined for nice functions $h\in
L^p_{\AB}(\Omega)$, for instance, for all compactly supported
continuous functions from $\R_+$ to $L^p_{\B}(\Omega)$. By the
density of all such functions in $L^p_{\AB}(\Omega)$, to prove the
boundedness of $Q$ we need only to estimate the relative norm of
$Qh$ for all such $h$.

\begin{theorem}\label{projection semigroup}
Let $\B, p, q$ be as in Theorem \ref{type-cotype semigroup}. Then
for any (nice) $h$
 $$\|{\mathfrak G}_q(Qh)\|_{L^p(\Omega)}\leq C_{p, q}\;
 \|h\|_{L^p_{\AB}(\Omega)}\ .$$
Consequently, ${\mathfrak G}_qQ$ extends to a bounded operator
from $L^p_{\AB}(\Omega)$ to $L^p(\Omega)$ with norm controlled by
a constant depending only on $p$ and $q$.
\end{theorem}

Admitting this theorem, we can easily prove Theorem
\ref{type-cotype semigroup}.

\medskip

\noindent{\sc Proof of Theorem \ref{type-cotype semigroup}}. {\sf
i)} $\Rightarrow$ {\sf ii)}. The proof for this is similar to that
for  {\sf i)} $\Rightarrow$ {\sf ii)} in Theorem \ref{lusin type
semigroups}.

{\sf ii)} $\Rightarrow$ {\sf i)}.  Fix an $f\in L^p_{\B}(\Omega)$.
Choose $h\in L_{L^{q'}_{\B^*}(\R_+,\frac{dt}{t})}^{p'}(\Omega)$ of
unit norm such
 $$ \|{\mathfrak G}_qf\|_{L^p(\Omega)}
 =\int_\Omega\int_{\R_+}\langle t\frac{\pa \Pt_tf}{\pa t}(x),\;
  h_t(x)\rangle\,\frac{dt}{t}\,dx.$$
Now we apply Theorem \ref{projection semigroup} to $\B^*, p'$ and
$q'$. (We may assume $f$ and $h$ are nice enough to legitimate the
calculations below.) We have, by hypothesis {\sf ii)} and Theorem
\ref{projection semigroup}
\begin{eqnarray*}
 \|{\mathfrak G}_qf\|_{L^p(\Omega)}
 &=&\int_{\R_+}\int_\Omega\langle f(x),\;
 t\frac{\pa \Pt_th_t}{\pa t}(x)\rangle\,\frac{dt}{t}\,dx,\\
 &=&\int_\Omega\langle f(x),\;
  Q(h)(x)\rangle\,dx,\\
 &\le&\|f\|_{L^p_{\B}(\Omega)}\|Q(h)\|_{L^{p'}_{\B^*}(\Omega)}\\
 &\le& C\|f\|_{L^p_{\B}(\Omega)}\|{\mathfrak
 G}_{q'}(Qh)\|_{L^{p'}(\Omega)}
 \le C'\|f\|_{L^p_{\B}(\Omega)}.
\end{eqnarray*}
This yields  {\sf i)}. \cqd

\medskip

As for  Theorem \ref{theorem on the means}, we shall reduce
Theorem \ref{projection semigroup} to an analogous inequality for
martingales via Rota's theorem. Let $\{E_n\}$ be a monotone
sequence of conditional expectations as in Lemma \ref{cesaro}. Let
us maintain the notations in that lemma and its proof. In the
remainder of this section $l^q$ denotes the usual $\ell^q$ space
over $\natt$ but with weight $\{\frac{1}{n}\}$, i.e., the norm of
a sequence $a$ is given by
 $$\|a\|_{l^q}=\Big(\sum_{n\ge 1}|a_n|^q\;\frac{1}{n}\Big)^{1/q}\
 .$$
The corresponding $\B$-valued version is $l_\B^q\ $, denoted by
$\D=l_\B^q$ in the sequel. Now we consider the discrete version of
$Q$ defined by (\ref{Q}). As before, the elements in $L^p_\D(M)$
are regarded as sequences with values in $L^p_\B(M)$. Given $h\in
L^p_\D(M)$ we define
 \begin{equation}\label{R}
 Rh=\sum_{n\ge 1} n\Delta\sigma_n(h_n)\;\frac{1}{n}\ ,
 \end{equation}
where $\Delta\sigma_n=\sigma_n - \sigma_{n-1}$. Recall that
$\displaystyle\sigma_n=\frac{E_0+\cdots+E_n}{n+1}$. $Rh$ is
clearly well-defined for finite sequences $(h_n)_n$.

\begin{lemma}\label{discrete projection}
Let $\B, p, q$ be as in Theorem \ref{type-cotype semigroup}. Let
$\{E_n\}$ be an arbitrary monotone sequence of conditional
expectations on a measure space $(M,dm)$. Then for any finite
sequence $h=(h_n)\in L^p_\D(M)$
 $$\big\|\big(n\Delta\sigma_n(Rh)\big)_{n\ge 1}\big\|_{L_\D^p(M)}
 \leq C_{p, q}\;\|h\|_{L^p_{\D}(M)}\ .$$
Consequently, $h\mapsto \big(n\Delta\sigma_n(Rh)\big)_{n\ge 1}$
extends to a bounded operator on $L^p_{\D}(M)$ with norm majorized
by a constant depending only on $p$ and $q$.
\end{lemma}

\noindent{\sc Proof}. Without loss of generality, we may assume
$\{E_n\}$ increasing. What we have to prove is the following inequality
 \begin{equation}\label{D}
 \Big\|\big[\sum_{m\ge1}\big\|\sum_{n\ge1}
 m\Delta\sigma_m\,\Delta\sigma_nh_n\big\|^q\frac{1}{m}
 \big]^{\frac{1}{q}}\Big\|_{L^p(M)}
 \le C_{p,q}\Big\|\big[\sum_{n\ge1}\|h_n\|^q\,\frac{1}{n}
 \big]^{\frac{1}{q}}\Big\|_{L^p(M)}.
 \end{equation}
 Given $m, n\ge 0$ we have
 \begin{equation}\label{ss}
 \Delta\sigma_m\,\Delta\sigma_n=\sigma_m\,\sigma_n -
 \sigma_{m-1}\,\sigma_n -\sigma_m\,\sigma_{n-1}
 +\sigma_{m-1}\,\sigma_{n-1}.
 \end{equation}
A simple calculation yields (with $m\wedge n=\min(m,n)$)
 $$\sigma_m\,\sigma_n=\frac{1}{(m+1)(n+1)}\;\sum_{i=0}^m\sum_{j=0}^n
 E_iE_j=\frac{1}{(m+1)(n+1)}\;\sum_{i=0}^{m\wedge
 n}(m+n-2i+1)E_i.$$
Setting $d_i=E_i-E_{i-1}$ (with the convention that $E_{-1}=0$)
and using Abel summation, we get
 $$\sigma_m\,\sigma_n=\sum_{j=0}^{m\wedge
 n}\Big(1-\frac{j(m+n+2-j)}{(m+1)(n+1)}\Big)d_j\ .$$
Thus by (\ref{ss})
 \begin{equation}\label{difference}
 \Delta\sigma_m\,\Delta\sigma_n=\frac{1}{mn(m+1)(n+1)}
 \sum_{j=1}^{m\wedge n}j^2\,d_j\ .
 \end{equation}
To prove \eqref{D} we use martingale transforms with vector-valued kernel
as in \cite{ma-to:00}. See the end of this paper (after Theorem
\ref{martingale th}) for a brief discussion on this subject. The
martingales we consider here are those defined on $(M,dm)$
relative to $\{E_n\}$ with values in $\D=l^q_\B$. Because of
(\ref{difference}), we want to express the mapping
 $$T: h=\{h_n\}_{n\ge1}\mapsto\Big\{\sum_{n\ge1}
 \frac{1}{n(m+1)(n+1)}\sum_{j=1}^{m\wedge n}
 j^2\,d_j\, h_n\Big\}_{m\ge1}$$
as a martingale transform, namely, we have to find a multiplying
sequence $\{v_j\}$ such that
 $$Th=\sum_jv_j\,d_jh.$$
This is clear from the above definition of $Th$. In fact, each
$v_j$ is a constant (an element in ${\mathcal L}(\D)$) and its matrix
is given by
 $$A_j=\Big[\frac{j^2}{n(m+1)(n+1)}\ind_{[j,\;\infty)}(m)
 \,\ind_{[j,\;\infty)}(n)
 \Big]_{m\ge1, n\ge 1}\ .$$
More precisely, $v_j$ is the operator $A_j\otimes \mbox{Id}_\B$,
where $A_j$ is considered as an operator on $l^q$. Therefore,
 $$\|v_j\|_{{\mathcal L}(\D)}=j^2\,
 \big(\sum_{m\ge j}\frac{1}{m(m+1)^q}\big)^{1/q}
 \big(\sum_{n\ge j}\frac{1}{n(n+1)^{q'}}\big)^{1/q'}
 \le C_0\ .$$
Therefore by \cite {ma-to:00}, it suffices to prove (\ref{D}) for
$p=q$. This is the main part of the proof. In the following $C>0$
denotes a constant which may depend on $q$  and vary from line to
line.

Let us first rewrite (\ref{D}) in the case $p=q$ (by using
(\ref{difference})):
 \begin{equation}\label{Dbis}
 L\equiv\int_M\sum_{m\ge1}\big\|\sum_{n\ge1}\frac{1}{n(m+1)(n+1)}
 \sum_{j=1}^{m\wedge n}j^2\,d_j\, h_n\big\|^q\,\frac{1}{m}
 \le C\int_M\sum_{n\ge1}\|h_n\|^q\,\frac{1}{n}.
 \end{equation}
It is easy to see that it is enough to prove the above inequality (\ref{Dbis})
for $m\wedge n\geq 5$. We have
 \begin{eqnarray}\label{L}
 \begin{array}{ccl}
 L&\le&\begin{displaystyle}
 \int_M \sum_{m\ge 5}\frac{1}{m^{q+1}}\,\big\|\sum_{5\leq n\le m}
 \frac{1}{n(n+1)}\sum_{j=1}^{n}j^2\,d_j\, h_n\big\|^q
 \end{displaystyle}\\
 & &+\begin{displaystyle}
 \;\int_M \sum_{m\ge5}\frac{1}{m^{q+1}}\,\big\|\sum_{n>m}
 \frac{1}{n(n+1)}\sum_{j=1}^{m}j^2\,d_j\, h_n\big\|^q
 \equiv A+B.
\end{displaystyle}
 \end{array}
 \end{eqnarray}
Let us first estimate $A$. To this end, we use the notations
$\Delta_k$ and $J_n$ introduced during the proof of Lemma
\ref{cesaro}. Note that for $2^{k-1}<j\le 2^{k}$
 $$d_j=E_j(\Delta_k)-E_{j-1}(\Delta_k).$$
Thus by Abel's summation we have
 $$\sum_{j=2^{k-1}+1}^{2^{k}}j^2\,d_j
 =2^{2k}\Delta_k-
 \sum_{j=2^{k-1}+1}^{2^{k}-1}(2j+1)E_j(\Delta_k),\quad k\geq 2$$
and so for $n\geq 5$
 \begin{eqnarray}\label{D1}
 \begin{array}{ccl}
 \sum_{j=1}^{n}j^2\,d_j
 &=&
 \begin{displaystyle}d_1+4d_2+\sum_{k=2}^{J_n}\sum_{j=2^{k-1}+1}^{2^{k}}j^2\,d_j
 +\sum_{j=2^{J_n}+1}^{n}j^2\,d_j
 \end{displaystyle}\\
 &=&\begin{displaystyle}
 d_1+4d_2+
 \sum_{k=2}^{J_n}2^{2k}\Delta_k
 -\sum_{k=2}^{J_n}\sum_{j=2^{k-1}+1}^{2^{k}-1}(2j+1)E_j(\Delta_k)
 \end{displaystyle}\\
 &&+\begin{displaystyle}\;n^2E_n(\Delta_{J_n+1})
 -\sum_{j=2^{J_n}+1}^{n-1}(2j+1)E_j(\Delta_{J_n+1}).
 \end{displaystyle}
 \end{array}
 \end{eqnarray}
Inserting this decomposition of $\sum_{j=1}^{n}j^2\,d_j$ in the
expression of $A$ in (\ref{L}) and by triangle inequality, we see
that $A$ is majorized by a sum of five terms, $A_0+A_1+A_2+A_3+A_4$,
corresponding respectively to the five terms of the last member in
(\ref{D1}). The term $A_0$ is handled easily.
Fix $0<\alpha<1/2$ and put $\beta=1-\alpha$. For $A_1$
we  have
\begin{eqnarray*}
 A_1
 &=&\int_M \sum_{m\ge5}\frac{1}{m^{q+1}}\,\big\|\sum_{5\leq n\le m}
 \frac{1}{n(n+1)}\sum_{k=2}^{J_n}2^{2k}\Delta_k\, h_n\big\|^q\\
 &=&\int_M \sum_{m\ge5}\frac{1}{m^{q+1}}\,\big\|\sum_{k=2}^{J_m}
 2^{2k}\Delta_k\big[\sum_{n:\;J_n\ge k}\frac{1}{n(n+1)}
 \,h_n\big]\big\|^q\\
 &\le&C\sum_{m\ge5}\frac{m^{2\alpha q}}{m^{q+1}}\sum_{k=2}^{J_m}
 2^{2\beta qk}\int_M \,\big\|\Delta_k\big[\sum_{n:\;J_n\ge k}\frac{1}{n(n+1)}
 \,h_n\big]\big\|^q\\
 &\le&C\int_M \sum_{m\ge5}\frac{m^{2\alpha q}}{m^{q+1}}\sum_{k=2}^{J_m}
 2^{2\beta qk}\,\big\|\sum_{n:\;J_n\ge k}\frac{1}{n(n+1)}
 \,h_n\big\|^q\\
 &\le&C\int_M \sum_{m\ge5}\frac{1}{m^{(1-2\alpha)q+1}}\sum_{k=2}^{J_m}
 2^{2\beta qk}\,2^{-k}\sum_{n:\;J_n\ge k}\frac{1}{n^q}
 \|h_n\|^q\\
 &\le&C\int_M \sum_{m\ge5}\frac{1}{m^{(1-2\alpha)q+1}}\sum_{5\leq n\le m}
 \frac{1}{n^q}\|h_n\|^q\sum_{k=2}^{J_n}2^{(2\beta q-1)k}\\
 &\le&C\int_M \sum_{n\ge 5}n^{2\beta q -q-1}\|h_n\|^q
 \sum_{m\ge n}\frac{1}{m^{(1-2\alpha)q+1}} \ \ (\mbox{since }
 2\beta q-1>0)\\
 &\le&C\int_M \sum_{n\ge 1}\frac{\|h_n\|^q}{n} \ \  (\mbox{since }
 1-2\alpha>0).
\end{eqnarray*}
We pass to $A_2$:
\begin{eqnarray*}
 A_2
 &=&\int_M \sum_{m\ge5}\frac{1}{m^{q+1}}\,\big\|\sum_{5\leq n\leq m}
 \frac{1}{n(n+1)}\sum_{k=2}^{J_n}\sum_{j=2^{k-1}+1}^{2^{k}-1}
 (2j+1)E_j\Delta_k\,h_n \big\|^q\\
 &=&\int_M \sum_{m\ge5}\frac{1}{m^{q+1}}\,\big\|\sum_{k=2}^{J_m}
 \big[\sum_{j=2^{k-1}+1}^{2^{k}-1}
 (2j+1)E_j\Delta_k\big]\big(\sum_{n:\;J_n\ge k}\frac{1}{n(n+1)}
 \,h_n\big)\big\|^q\\
\end{eqnarray*}
Let us handle the expression concerning the internal norm. As for
$A_1$ previously, we have
\begin{eqnarray*}
 &&\big\|\sum_{k=2}^{J_m}
 \big[\sum_{j=2^{k-1}+1}^{2^{k}-1}
 (2j+1)E_j\Delta_k\big]\big(\sum_{n:\;J_n\ge k}\frac{1}{n(n+1)}
 \,h_n\big)\big\|^q\\
 &\le& C\; m^{2\alpha q}\sum_{k=2}^{J_m}2^{2\beta qk}
 \big\|\sum_{j=2^{k-1}+1}^{2^{k}-1}2^{-2k}(2j+1) E_j\Delta_k
 \big(\sum_{n:\;J_n\ge k}\frac{1}{n(n+1)}\,h_n\big)\big\|^q\\
 &\le& C\; m^{2\alpha q}\sum_{k=2}^{J_m}2^{2\beta qk}
 \sum_{j=2^{k-1}+1}^{2^{k}-1}2^{-2k}(2j+1)\big\|E_j\Delta_k
 \big(\sum_{n:\;J_n\ge k}\frac{1}{n(n+1)}\,h_n\big)\big\|^q\ .
\end{eqnarray*}
Therefore
\begin{eqnarray*}
 &&\int_M\big\|\sum_{k=2}^{J_m}
 \big[\sum_{j=2^{k-1}+1}^{2^{k}-1}
 (2j+1)E_j\Delta_k\big]\big(\sum_{n:\;J_n\ge k}\frac{1}{n(n+1)}
 \,h_n\big)\big\|^q\\
 &\le& C\; m^{2\alpha q}\sum_{k=2}^{J_m}2^{2\beta qk}
 \sum_{j=2^{k-1}+1}^{2^{k}-1}2^{-2k}(2j+1)\int_M\big\|E_j\Delta_k
 \big(\sum_{n:\;J_n\ge k}\frac{1}{n(n+1)}\,h_n\big)\big\|^q\\
 &\le& C\; m^{2\alpha q}\sum_{k=2}^{J_m}2^{2\beta qk}
 \sum_{j=2^{k-1}+1}^{2^{k}-1}2^{-2k}(2j+1)\int_M\big\|
 \sum_{n:\;J_n\ge k}\frac{1}{n(n+1)}\,h_n\big\|^q\\
 &\le& C \int_Mm^{2\alpha q}\sum_{k=2}^{J_m}2^{2\beta qk}
 \big\|\sum_{n:\;J_n\ge k}\frac{1}{n(n+1)}\,h_n\big\|^q\ .
\end{eqnarray*}
Combining the preceding inequalities, we get
\begin{eqnarray*}
 A_2
 &\le&C\int_M \sum_{m\ge5}\frac{1}{m^{(1-2\alpha)q+1}}\sum_{k=2}^{J_m}
 2^{2\beta qk}\,\big\|\sum_{n:\;J_n\ge k}\frac{1}{n(n+1)}
 \,h_n\big\|^q
 \le C\int_M \sum_{n\ge 1}\frac{\|h_n\|^q}{n}\ .
\end{eqnarray*}
It is easier to estimate $A_3$ and $A_4$. Indeed, for $A_3$ we
have
\begin{eqnarray*}
 A_3
 &=&\int_M \sum_{m\ge5}\frac{1}{m^{q+1}}\,\big\|\sum_{5\leq n\le m}
 \frac{1}{n(n+1)}\,n^{2}E_n\Delta_{J_n+1}\, h_n\big\|^q\\
 &\le&C \sum_{m\ge5}\frac{1}{m^{2}}\,\sum_{5\leq n\le m}
 \int_M \|E_n\Delta_{J_n+1}\, h_n\big\|^q\le C\int_M \sum_{n\ge 1}\frac{\|h_n\|^q}{n} .
\end{eqnarray*}
$A_4$ is similarly estimated. Therefore, we get
\begin{equation}\label{A}
 A\le C\|h\|_{L^q_\D(M)}^q\ .
\end{equation}

Now we turn to $B$. This time we fix $\alpha$ such that
$\frac{q-1}{2q}<\alpha<\frac{q-1}{q}$ and put again $\beta
=1-\alpha$. As for $A$ previously, using (\ref{D1}) with $n$
replaced by $m$, we see that $B$ is less than or equal to
$B_0+B_1+B_2+B_3+B_4$, corresponding to the decomposition in
(\ref{D1}). Thus
\begin{eqnarray*}
 B_1
 &=&\int_M \sum_{m\ge5}\frac{1}{m^{q+1}}\,\big\|\sum_{n> m}
 \frac{1}{n(n+1)}\sum_{k=2}^{J_m}2^{2k}\Delta_k\, h_n\big\|^q\\
 &=&\int_M \sum_{m\ge5}\frac{1}{m^{q+1}}\,\big\|\sum_{k=2}^{J_m}
 2^{2k}\Delta_k\big[\sum_{n>m}\frac{1}{n(n+1)}
 \,h_n\big]\big\|^q\\
 &\le&C\sum_{m\ge5}m^{q-3}\sum_{k=2}^{J_m}
 2^{2k}\int_M \big\|\Delta_k\big[\sum_{n>m}\frac{1}{n(n+1)}
 \,h_n\big]\big\|^q\\
 &\le&C\int_M \sum_{m\ge5}m^{q-1}
 \big(\sum_{n>m}\frac{1}{n^{2\alpha q'}}\big)^{q-1}
 \sum_{n>m}\frac{1}{n^{2\beta q}}\|h_n\big\|^q\\
 &\le&C\int_M \sum_{m\ge5}m^{2(q-1)-2\alpha q}\sum_{n>m}
 \frac{1}{n^{2\beta q}}\|h_n\big\|^q\ \ (\mbox{since }2\alpha q'>1)\\
 &\le&C\int_M \sum_{n\ge 1}\frac{\|h_n\|^q}{n} \ \  (\mbox{since }
 q-1>\alpha q).
\end{eqnarray*}
In a similar way (also as for $A_2, A_3$ and $A_4$), we obtain the
same bound for $B_2, B_3$ and $B_4$. Therefore, we have
\begin{equation}\label{B}
 B\le C\|h\|_{L^q_\D(M)}^q\ .
\end{equation}
Finally combining (\ref{L}), (\ref{A}) and (\ref{B}), we get
(\ref{Dbis}). Therefore, the lemma is proved.\cqd

\medskip

\noindent{\sc Proof of Theorem \ref{projection semigroup}.} For
notational simplicity we  put
 $$\Phi_t=t\frac{\pa P_t}{\pa t}\quad\mbox{and}\quad
 \Psi_t=t\frac{\pa M_t}{\pa t}.$$
(Recall that $M_t=\frac{1}{t}\int_0^tT_s\,ds$.) By Lemma
\ref{poisson-mean},
 $$\int_\Omega\big[\int_0^\infty\big\|\int_0^\infty\Phi_s
 \Phi_th_t\,\frac{dt}{t}\big\|^q\big]^{\frac{p}{q}}\,\frac{ds}{s}
 \le C^p\int_\Omega\big[\int_0^\infty\big\|\int_0^\infty\Psi_s
 \Phi_th_t\,\frac{dt}{t}\big\|^q\big]^{\frac{p}{q}}\,\frac{ds}{s}.$$
By (\ref{poisson-mean1}),
\begin{eqnarray*}
 \int_0^\infty\Psi_s\Phi_th_t\,\frac{dt}{t}
  &=&-2\int_0^\infty\int_0^\infty\big(\frac{u}{t^2}\big)^2
  \varphi'(\frac{u}{t^2})\Psi_s\Psi_uh_t\,\frac{du}{u}\,\frac{dt}{t}\\
  &=&\int_0^\infty\int_0^\infty\big[t^2
  \varphi'(t)\big]\Psi_s\Psi_uh_{u^{1/2}t^{-1/2}}
  \,\frac{du}{u}\,\frac{dt}{t}\ .
\end{eqnarray*}
By the triangle inequality (recalling that
$\AB=L^q_\B(\R_+,\frac{ds}{s})$),
\begin{eqnarray*}
 \Big\|\int_0^\infty\Psi_\cdot\,\Phi_th_t\,\frac{dt}{t}\Big\|_
 {L^p_{\AB}(\Omega)}
  &\le&\int_0^\infty\Big\|\int_0^\infty\Psi_\cdot\,
  \Psi_u\,h_{u^{1/2}t^{-1/2}}
  \,\frac{du}{u}\Big\|_{L^p_{\AB}(\Omega)}
  t^2|\varphi'(t)|\frac{dt}{t}\ .
\end{eqnarray*}
Now using the same discretization arguments as in the proof of
Theorem \ref{theorem on the means}, we deduce from Lemma
\ref{discrete projection} that
 \begin{eqnarray*}
 \Big\|\int_0^\infty\Psi_\cdot\,\Psi_u\,h_{u^{1/2}t^{-1/2}}
  \,\frac{du}{u}\Big\|^p_{L^p_{\AB}(\Omega)}
 &\le& C^p \int_\Omega\big[\int_0^\infty\|h_{u^{1/2}t^{-1/2}}\|^q
 \frac{du}{u}\big]^{\frac{p}{q}}\\
 &=& 2^{p/q}C^p \|h\|^p_{L^p_\AB(\Omega)}\ .
\end{eqnarray*}
Combining the preceding inequalities, we obtain
$$
 \Big\|\int_0^\infty\Phi_\cdot\,\Phi_th_t\,\frac{dt}{t}\Big\|_
 {L^p_{\AB}(\Omega)}
 \le C'\|h\|_{L^p_\AB(\Omega)}\ .
$$
This is the desired inequality. Thus we have achieved the proof of
Theorem \ref{projection semigroup}.\cqd


\section{Poisson semigroup on $\R^n$}

This section and the next are devoted to the study of the
Littlewood-Paley $g$-function on $\R^n$ in the vector-valued case.
Our main goal is to prove the implications {\sf ii)} $\Rightarrow$
{\sf i)} in Theorems \ref{lusin semigroups} and \ref{lusin type
semigroups} in the particular case of the Poisson semigroup on
$\R^n$. This section collects some results on the $g$-function
operator on $\R^n$, represented as a Calder\'on-Zygmund operator.
It can be considered as preparatory, although some of these
results are of general interest. The proof of the mentioned
implications will be done in the next section.

Let $\B$ be a Banach space and $1<q<\infty$. Recall the
generalized  Littlewood-Paley $g$-function on $\R^n$:
$$
\g_q(f)(x)=\Big(\int_0^\infty t^q\|\nabla P_t*
f(x)\|_{\ell^2_\B}^q\,\frac{dt}{t}\Big)^{1/q},\quad x\in\R^n\ .
$$
It is often easier to consider the corresponding $g$-function
defined only by the derivative in time, which is the following
$$
\g^1_q(f)(x)=\Big(\int_0^\infty t^q\big\|\frac{\partial
P_t}{\partial t} *f(x)\big\|_{\B}^q\,\frac{dt}{t}\Big)^{1/q}\ .
$$
Similarly, we define $\g_q^2(f)$ as the part of $\g_q(f)$
corresponding to the gradient in the space variable:
$$
\g^2_q(f)(x)=\Big(\int_0^\infty t^q\|\nabla_x P_t*
f(x)\|_{\ell^2_\B}^q\,\frac{dt}{t}\Big)^{1/q},
$$
where
$$\nabla_x=\big(\frac{\pa}{\pa x_1},\;...\;, \frac{\pa}{\pa
x_n}\big).$$
 These $g$-functions can be treated as
Calder{\'o}n-Zygmund operators. To this end we first recall
briefly the definition of these operators. Given a pair of Banach
spaces $\B_1$ and $\B_2$, a linear operator $T$ is a
Calder\'on-Zygmund operator on $\R^n$, with associated
Calder\'on-Zygmund kernel $k$ if $T$ maps
$L^\infty_{c,\B_1}(\R^n)$, the space of the essentially bounded
functions on $\R^n$ with compact support,  into $\B_2$-valued
strongly measurable functions on $\R^n$, and for any function
$f\in L^\infty_{c,\B_1}(\R^n)$ we have
$$ Tf(x) = \int_{\R^n}k(x,y)f(y)dy, \quad\mbox{for a. e. }x \mbox{
outside the support of } f,$$
where the kernel, $k(x,y) \in {\mathcal L}({\B_1},{\B_2})$ satisfies
\begin{itemize}
\item[{\sf a)}] $\|k(x,y)\| \leq C|x-y|^{-n}$
\item[{\sf b)}] $\|\nabla_xk(x,y)\| + \|\nabla_yk(x,y)\|
\leq C|x-y|^{-(n+1)}$
\end{itemize}
We shall always assume that there is $\Lambda\in {\mathcal
L}({\B_1},{\B_2})$ such that $T(c)\equiv\Lambda(c)$ for all $c\in
\B_1$.

\medskip

Let us recall the $\BMO$ and $H^1$ spaces on $\R^n$. Let $\B$ be a
Banach space. $\BMO_\B(\R^n)$ is the space of $\B$-valued
functions $f$ defined on $\R^n$ such that
$$
\|f\|_{\BMO_\B(\R^n)}
=\sup_Q\frac{1}{|Q|}\int_Q\|f(x)-f_Q\|_\B\,dx<\infty,
$$
where $f_Q=\frac{1}{|Q|}\int_Qf(x)\,dx$ and the supremum is taken
over the cubes $Q\subset\R^n$. The space $H^1$ is defined in the
atomic sense. Namely, we say that a function $a\in
L^\infty_\B(\R^n)$ is a $\B$-{\em atom} if there exists a cube
$Q\subset\R^n$ containing the support of $a$, and such that
$\|a\|_{L^\infty_\B(\R^n)}\leq |Q|^{-1}$, and $\int_Qa(x)\,dx=0$.
Then, we say that a function $f$ is in $H^1_\B(\R^n)$ if it admits
a decomposition $f=\sum_i\lambda_ia_i$, where $a_i$ are
$\B$-valued atoms and $\sum_i|\lambda_i|<\infty$. We define
$\|f\|_{H^1_\B}=\inf\big\{\sum_i|\lambda_i|\big\}$, where the
infimum runs over all those such decompositions (see
\cite{bl-ga:87}).

\medskip

The following theorem is a kind of folklore. We give a sketch of
its proof for the convenience of the reader. $\BMO_{c,\B}(\R^n)$
denotes the subspace of $\BMO_{\B}(\R^n)$ consisting of
functions with compact support.

\begin{theorem}\label{singular integral}
Let $\B_1, \; \B_2$ be two Banach spaces and $T$ a
Calder\'on-Zygmund operator with an associated kernel $k$ as
above. Let $S$ be defined as $S(f) = \|T(f)\|_{\B_2}$. Then, the
following statements are equivalent
\begin{itemize}
\item[{\sf i)}] $T$ maps $L^\infty_{c,\B_1}(\R^n)$ into $BMO_{\B_2}(\R^n)$.
\item[{\sf ii)}] $S$  maps $L^\infty_{c,\B_1}(\R^n)$ into $\BMO(\R^n)$
\item[{\sf iii)}] $T$  maps $H^1_{\B_1}(\R^n)$ into $L^1_{\B_2}(\R^n)$.
\item[{\sf iv)}] $T$ maps $L^p_{\B_1}(\R^n)$ into $L^p_{\B_2}(\R^n)$ for any
(or equivalently, for some) $p\in (1,\infty)$.
\item[{\sf v)}] $T$ maps $L^1_{\B_1}(\R^n)$ into weak-$L^{1}_{\B_2}(\R^n)$.
\item[{\sf vi)}] $T$ maps $\BMO_{c,\B_1}(\R^n)$ into
$\BMO_{\B_2}(\R^n)$.
\item[{\sf vii)}] $S$ maps $\BMO_{c,\B_1}(\R^n)$ into $\BMO(\R^n)$.
\end{itemize}
\end{theorem}

\noindent {\sc Proof}. The structure of the proof is the
following: first,
 {\sf i)} $\Rightarrow$ {\sf ii)} $\Rightarrow$ {\sf iii)} $\Rightarrow$ {\sf i)}.
 Then, we prove {\sf i)} $\Rightarrow$ {\sf iv)}
 $\Rightarrow$ {\sf vi)} $\Rightarrow$ {\sf vii)} $\Rightarrow$ {\sf ii)}
and    {\sf iv)} $\Rightarrow$ {\sf v)} $\Rightarrow$ {\sf i)}.

The fact that $L^\infty_{\B_1}$ is contained in $\BMO_{\B_1}$
gives that {\sf vii)} implies {\sf ii)}. Since the norm of  a
function  in $\BMO_\B$  is a function in $\BMO$, then we have {\sf
i)} $\Rightarrow$ {\sf ii)} and {\sf vi)} $\Rightarrow$ {\sf vii)}.

To get {\sf ii)} $\Rightarrow$ {\sf iii)} and {\sf iii)}
$\Rightarrow$ {\sf i)}, we can proceed as in \cite [p.49] {jou}
with minor modifications due to considering the operator
$Sf=\|Tf\|_{\B_2}$.

 As we already know that {\sf i)}$\Rightarrow$ {\sf iii)},
we can apply  interpolation (see \cite{bx}) and we have that $T$
maps $L^p_{\B_1}$ into $L^p_{\B_2}$ for $1<p< \infty$, so we have
{\sf i)} $\Rightarrow$ {\sf iv)}.

The proof of {\sf iv)} $\Rightarrow$ {\sf vi)} is where the
condition that $T(c)(x)=\Lambda(c)$ plays a role. Let $f$ be a
function in $BMO_{c,\B_1}$. Given a cube $Q$ with center $x_0$,
let $\tilde Q$ be its doubled cube. We decompose
\begin{eqnarray*}
\lefteqn{\frac1{|Q|}\int_Q \|Tf(x)-Tg_2(x_0)\|_{\B_2}\,dx\leq}\\
&\leq& \frac1{|Q|}\int_Q \|Tg_1(x)\|_{\B_2}\,dx+ \frac1{|Q|}\int_Q
\|Tg_2(x)-Tg_2(x_0)\|_{\B_2}\,dx,
\end{eqnarray*}
where $ f = g_1+g_2$, $g_1= (f-f_Q)\ind_{\tilde{Q}}$ and
$g_2=(f-f_Q)\ind_{\R^n\setminus \tilde{Q}} +f_Q$. By using Jensen
and the $L^p$ boundedness of $T$, we have
\begin{eqnarray*}
\frac{1}{|Q|}\int_{Q}\|Tg_1(x)\|\,dx &\leq&
\Big(\frac{1}{|Q|}\int_{Q}\|Tg_1(x)\|^p\,dx\Big)^{1/p} \leq C
\Big(\frac{1}{|Q|}\int_{\R^n}\|g_1(x)\|_{\B_1}^p\,dx\Big)^{1/p} \\
&=& C\Big(\frac{1}{|Q|}\int_{\tilde
Q}\|f(x)-f_Q\|_{\B_1}^p\,dx\Big)^{1/p} \leq
C\|f\|_{\BMO_{\B_1}},
\end{eqnarray*}
where in the last inequality we have used the John-Niremberg
theorem. On the other hand, using  $T(c)(x) = \Lambda(c)$, we have
\begin{eqnarray*}
Tg_2(x)-Tg_2(x_0) &=&
T( (f-f_Q)\ind_{\R^n\setminus\tilde{Q}})(x)
-T( (f-f_Q)\ind_{\R^n\setminus\tilde{Q}})(x_0) \\
&=& \int_{\R^n\setminus\tilde{Q}}(k(x,y)-k(x_0,y))(f(y)-f_Q)\,dy.
\end{eqnarray*}
Now, using the hypothesis on the kernel $k$ we have that for $x\in
Q$, $y\in\R^n\setminus\tilde Q$, $\|k(x,y)-k(x_0,y)\|\leq
\frac{|x-x_0|}{|y-x_0|^{n+1}}$, and therefore
\begin{eqnarray*}
\|Tg_2(x)-Tg_2(x_0)\|_{\B_2} &\le&
\sum_{j=1}^\infty 2^{-j}\int_{2^jQ\setminus 2^{j-1}Q}
\frac{1}{|y-x_0|^n}\|f(y)-f_Q\|_{\B_1}\,dy \\
&\leq& \sum_{j=1}^\infty
2^{-j}\frac1{|2^jQ|}\int_{2^jQ}\|f(y)-f_Q\|_{\B_1}\,dy \leq C
\|f\|_{\BMO_{\B_1}}.
\end{eqnarray*}

By \cite [Theorem V.3.4] {ga-rdf:85}, the strong type $(p,p)$
implies that $T$ is of weak type $(1,1)$. This gives
{\sf iv)} $\Rightarrow$ {\sf v)}.
>From {\sf v)}, statement {\sf i)} can be
achieved by using an slight modification of the argument given in
the proof of Lemma 5.11 of \cite[pg. 199]{ga-rdf:85}. The key is
using Kolmogorov's inequality \cite[Lemma V.2.8]{ga-rdf:85} relating
$L^{1,\infty}$ norm with $L^q$ norm for $0<q<1$ and the fact that
$BMO_q=BMO$.\cqd

\medskip

It is well known that the various Littlewood-Paley $g$-functions
can be expressed as Calder{\'o}n-Zygmund operators with regular
vector-valued  kernels (see \cite{xu} for the case of the torus;
also see \cite{st:2} for the scalar case). Therefore, we
immediately get the following

\begin{corollary}\label{littlewood-paley}
Given a Banach space $\B$, $q\in(1,\infty)$, the following
statements are equivalent.
\begin{itemize}
\item[{\sf i)}] $\g_q$ maps $L^\infty_{c,\B}(\R^n)$ into $\BMO(\R^n)$.
\item[{\sf ii)}] $\g_q$ maps $H^1_\B(\R^n)$ into $L^1(\R^n)$.
\item[{\sf iii)}] $\g_q$ maps $L^p_{\B}(\R^n)$ into $L^p(\R^n)$
for any (equivalently for some) $p\in(1,\infty)$.
\item[{\sf iv)}] $\g_q$ maps $\BMO_{c,\B}(\R^n)$ into $\BMO(\R^n)$.
\item[{\sf v)}] $\g_q$ maps $L^{1}_\B(\R^n)$ into $L^{1,\infty}(\R^n)$.
\end{itemize}
These statements are also equivalent if we replace $\g_q$ by $\g^1_q$ or $\g^2_q$.
\end{corollary}

The following result reduces the boundedness on $\g_q$, $\g^1_q$
and $\g^2_q$ to that on one of them.

\begin{proposition}\label{time-space}
Let $\B$ be a Banach space and $p,q\in(1,\infty)$. Then for any
$f\in L_\B^p(\R^n)$
$$\|\g^1_q(f)\|_{L^p(\R^n)}\approx \|\g^2_q(f)\|_{L^p(\R^n)}\ ,$$
where the equivalence constants depend only on $p, q$ and $n$.
\end{proposition}

\noindent {\sc Proof}.  Set, for simplicity,
$$
 \pa=\frac{\pa}{\pa t},\
 \pa_i=\frac{\pa}{\pa x_i}\quad\mbox{and}\quad
 \Phi_t=t\;\pa P_t,\ \Phi^i_t=t\;\pa_i P_t,\quad i=1,\; ...\;, n.
$$
Given $t=t_1+t_2$ we have
 $$P_t*f=P_{t_1}*P_{t_2}*f.$$
Differentiating the two sides first in $t_2$ and then in $x_i$, we
get
\begin{equation}\label{tx}
 t^2\,\pa_i\pa P_{2t}*f=\Phi^i_{t}*\Phi_{t}*f.
\end{equation}
(Here $\pa_i\pa P_{2t}*f=\big[\pa_{x_i}\pa_s
(P_{s}*f(x))\big]\big|_{s=2t}$.) We use again the singular
integral theory. Let $\AB=L^q_\B(\R_+,\;\frac{dt}{t})$. Given
$x\in\R^n$ let $k(x): \AB\to\AB $ be the operator defined by
$k(x)(\varphi)(t)=\Phi^i_{t}(x)\cdot \varphi(t)$ for
$\varphi\in\AB$. It is easy to check that $(x, y)\mapsto k(x-y)$
is a Calder\'on-Zygmund kernel (satisfying additionally the
condition in Theorem \ref{singular integral}). Let $T$ be the
associated operator. We claim that $T$ is bounded on
$L_{\AB}^q(\R^n)$. Indeed, fix $h\in L_{\AB}^q(\R^n)$. Note that
$h$ can be regarded as a function of two variables $(x, t)\in \R^n
\times\R_+$. Then
$$\big\|T(h)\big\|^q_{L_{\AB}^q(\R^n)}
 =\int_{\R^n}\int_{\R_+}\Big\|\int_{\R^n}\Phi_t^i(y)h(x-y,t)
 \,dy\Big\|^q\;\frac{dt}{t}\;dx.$$
However,
\begin{eqnarray*}
 \big(\int_{\R^n}\big\|\int_{\R^n}\Phi_t^i(y)h(x-y,t)
 \,dy\big\|^q\;dx\big)^{\frac{1}{q}}
 &\le&\int_{\R^n}\big|\Phi_t^i(y)\big|\big(\int_{\R^n}\|h(x-y,t)\|^q
 \,dx\big)^{\frac{1}{q}}dy\\
 &\le& C \big(\int_{\R^n}\|h(x,t)\|^q \,dx\big)^{\frac{1}{q}}\ .
\end{eqnarray*}
 Therefore,
$$\big\|T(h)\big\|_{L_{\AB}^q(\R^n)}\le C
 \|h\|_{L_{\AB}^q(\R^n)}\ .$$
This gives our claim. Then by Theorem \ref{singular integral}, we
deduce that $T$ is bounded on $L_{\AB}^p(\R^n)$ for all
$p\in(1,\infty)$.

Applying this boundedness of $T$ to $h(x, t)=\Phi_{t}*f(x)$ and
using (\ref{tx}), we get
$$
 \int_{\R^n}\big(\int_{\R_+}\|t^2\,\pa_i\pa P_{t}*f(x)\|^q\;
 \frac{dt}{t}\big)^{p/q}dx
 \le C^p\;\|\g^1_q(f)\|^p_{L^p(\R^n)}\ .$$
Assume, without loss of generality, that $f$ is good enough such
that $\pa_i P_{t}*f(x)\to 0$ as $t\to\infty$ (for instance, $f$ is
compactly supported). Then
 $$\pa_iP_{t}*f(x)=-\int_t^\infty \pa \pa_i P_{s}*f(x)ds.$$
Therefore
\begin{eqnarray*}
 \int_0^\infty\|t\pa_iP_{t}*f(x)\|^q\;\frac{dt}{t}
 &\le&\int_0^\infty t^q\big(\int_t^\infty s\|\pa\pa_i P_{s}*f(x)\|
 \;\frac{ds}{s}\big)^q\;\frac{dt}{t}\\
 &\le& C\int_0^\infty t^{q/2} \int_t^\infty s^{3q/2}
 \|\pa\pa_i P_{s}*f(x)\|^q \;\frac{ds}{s}\;\frac{dt}{t}\\
 &=&C'\int_0^\infty s^{2q}
 \|\pa\pa_i P_{s}*f(x)\|^q \;\frac{ds}{s}\ .
\end{eqnarray*} It then follows that
$$\int_{\R^n}\big(\int_{\R_+}\|t\,\pa_iP_{t}*f(x)\|^q\;
 \frac{dt}{t}\big)^{p/q}dx
 \le C^p\int_{\R^n}\big(\int_{\R_+}\|t^2\,\pa_i\pa P_{t}*f(x)\|^q\;
 \frac{dt}{t}\big)^{p/q}dx.$$
Therefore,
$$\int_{\R^n}\big(\int_{\R_+}\|t\,\pa_iP_{t}*f(x)\|^q\;
 \frac{dt}{t}\big)^{p/q}dx\le C^p\;\|\g^1_q(f)\|^p_{L^p(\R^n)}\ .$$
Adding the $n$ inequalities so obtained over $i=1,\;...\;,n$, we
get
 $$\|\g^2_q(f)\|_{L^p(\R^n)}
 \le C_{p,q,n}\;\|\g^1_q(f)\|_{L^p(\R^n)}\ .$$

To prove the reverse inequality, we first observe that the same
arguments as above show that
$$
 \int_{\R^n}\big(\int_{\R_+}\|t^2\,\pa_i^2 P_{t}*f(x)\|^q\;
 \frac{dt}{t}\big)^{p/q}dx
 \le C^p\int_{\R^n}\big(\int_{\R_+}\|t\,\pa_i P_{t}*f(x)\|^q\;
 \frac{dt}{t}\big)^{p/q}dx\ .$$
Then using the formula
 $$\pa^2 P_{t}*f(x)=-\sum_{i=1}^n\pa_i^2 P_{t}*f(x),$$
we get
$$
 \int_{\R^n}\big(\int_{\R_+}\|t^2\,\pa^2 P_{t}*f(x)\|^q\;
 \frac{dt}{t}\big)^{p/q}dx
 \le C^p\;\|\g^2_q(f)\|^p_{L^p(\R^n)}\ .$$
Finally, using the arguments in the second part of the above
proof, we can go down to $\pa P_{t}*f(x)$ to have
 $$\|\g^1_q(f)\|_{L^p(\R^n)}\le C_{p,q, n}\;
 \|\g^2_q(f)\|_{L^p(\R^n)}\ .$$
This completes the proof of the proposition.\cqd

\medskip

\noindent{\bf Remark.} The equivalence in Proposition
\ref{time-space} still holds when $L^p$ is replaced by
$L^{1,\infty}$.

\medskip

The last result of this section is a duality theorem for the
boundedness of the $g$-functions (Theorem \ref{type-cotype}
below). This theorem is a particular case of Theorem
\ref{type-cotype semigroup}. It is also the analogue for $\R^n$ of
\cite [Theorem 2.4] {xu} (for the torus). As in such a situation
the key is again the existence of a certain projection. For the
reader interested only in the $\R^n$ case, we include a proof for
this latter fact, which is much simpler than that of Theorem
\ref{projection semigroup}. Fix a Banach space $\B$ and
$q\in(1,\infty)$, and we keep the notation used in the previous
proof with $\AB=L_\B^q(\R_+,\; \frac{dt}{t})$. The projection in
question is defined as (recalling that $\Phi_t=t\pa P_t$)
 $$Q(h)(x)=\int_0^\infty \Phi_t*h(\cdot, t)(x)\frac{dt}{t}\ ,
 \quad x\in\R^n.$$
Note that $Q(h)$ is well defined for functions $h$ in a dense
family of $L^p_{\AB}(\R^n)$, for instance, for those which are
compactly supported continuous functions on $\R^n\times\R_+$.

The proof of the following lemma is an adaption for $\R^n$ of
\cite [Lemma 2.3] {xu}, so we are rather sketchy.

\begin{lemma} \label{projection}
Let $\B$ be a Banach space and $p,q\in(1,\infty)$. Let $Q$ be
defined as before. Then for any $\B$-valued continuous function
$h$ with compact support in $\R^n\times\R_+$
 $$\|\g^1_q(Q(h))\|_{L^p(\R^n)}\le
 C_{p,q} \|h\|_{L^p_{\AB}(\R^n)}\ .$$
Consequently, the map $h\mapsto \g^1_q(Q(h))$ extends to a bounded
map from $L^p_{\AB}(\R^n)$ to $L^p(\R^n)$.
\end{lemma}

\noindent {\sc Proof}. Set $f=Q(h)$. Then
\begin{eqnarray*}
 s\Phi_s*f
 &=&s\,\int_0^\infty \Phi_s*\Phi_t*h(\cdot, t)(x)\frac{dt}{t}\\
 &=&st\,\int_0^\infty \pa^2P_{s+t}*h(\cdot, t)(x)\frac{dt}{t}
 =\int_0^\infty k_{s,t}*h(\cdot, t)(x)\frac{dt}{t},
\end{eqnarray*}
where
 $$k_{s,t}=st\,\pa^2P_{s+t}=st\,\frac{\pa^2}{\pa
 u^2}P_{u}\Big|_{u=s+t}\ .$$
Now consider the operator $K(x): \AB\to\AB$ defined by
 $$K(x)(\varphi)(s)=\int_0^\infty
 k_{s,t}(x)\varphi(t)\frac{dt}{t}$$
for every $\varphi\in\AB$. Using the inequality
 $$|k_{s,t}(x)|\le C\,\frac{st}{(|x|+s+t)^{(n+2)}}$$
and a similar one for the derivative of $k_{s,t}(x)$ in $x$, one
can easily check that for any $x\in\R^n\setminus\{0\}$, $K(x)$ is
bounded and
 $$\|K(x)\|\le \frac{C}{|x|^n}\ ,\qquad
 \|\nabla K(x)\|\le \frac{C}{|x|^{n+1}}\ .$$
Thus $(x, y)\mapsto K(x-y)$ is a Calder\'on-Zygmund kernel. Hence
to prove the lemma it suffices to show that the singular integral
operator $h\mapsto K*h$ is bounded on $L_{\AB}^q(\R^n)$,  in
virtue of Theorem \ref{singular integral}. This is easily done as
follows. For $x\in\R^n$ and $s\in(0,\infty)$ we have
\begin{eqnarray*}
 \|K*h (x, s)\|
 &\le& \big(\int_0^\infty
 \int_{\R^n}|k_{s,t}(x-y)|\,dy\;\frac{dt}{t}\big)^{1/q'}\;\\
 &&\times\;\big(\int_0^\infty
 \int_{\R^n}|k_{s,t}(x-y)|\;\|h(y,t)\|^q\,dy\;\frac{dt}{t}\big)^{1/q}\\
 &\le& C \big(\int_0^\infty
 \int_{\R^n}|k_{s,t}(x-y)|\;\|h(y,t)\|^q\,dy\;\frac{dt}{t}\big)^{1/q}\
 .
\end{eqnarray*}
Therefore,
\begin{eqnarray*}
 \|K*h\|^q_{L_{\AB}^q(\R^n)}
 &\le& C^q \int_{\R^n}\int_0^\infty\int_0^\infty
 \int_{\R^n}|k_{s,t}(x-y)|\;\|h(y,t)\|^q\,dy\;
 \frac{dt}{t}\;\frac{ds}{s}\;dx\\
 &\le& C^q \int_{\R^n}\int_0^\infty [\int_{\R^n}\int_0^\infty
 \frac{st}{(|x-y|+s+t)^{(n+2)}}\frac{ds}{s} dx ]\|h(y,t)\|^q
 \frac{dt}{t} dy\\
 &\le& C_{q, n}^q\|h\|^q_{L_{\AB}^q(\R^n)}\ .
\end{eqnarray*}
This implies the desired boundedness of the singular integral on
$L_{\AB}^q(\R^n)$. \cqd

\medskip

\noindent{\bf Remark.} Lemma \ref{projection} holds as well for
$\g_q$ and $\g^2_q$ instead of $\g^1_q$. Moreover, the weak type
(1, 1) inequality is true too.

\medskip

>From Lemma \ref{projection} and using the arguments in the proof
of Theorem \ref{type-cotype semigroup}, we deduce the following

\begin{theorem}\label{type-cotype}
Let $\B$ be a Banach space and $q\in(1,\infty)$. Then the
following statements are equivalent:
\begin{itemize}
\item[{\sf i)}] One of the statements in Corollary
\ref{littlewood-paley} holds.
\item[{\sf ii)}] For every $p\in(1,\infty)$ (equivalently for some
$p\in(1,\infty)$) there is a constant $C$ such that
 $$\|f\|_{L^p_{\B^*}(\R^n)}\le C\,\|\g^1_{q'}\|_{L^p(\R^n)}\ ,
 \quad \forall\; f\in L^p_{\B^*}(\R^n).$$
\end{itemize}
\end{theorem}


\section{Poisson semigroup on $\R^n$ continued}

Our aim in this section is proving that in the definition of the
Lusin type or cotype the $G_q$-function on the torus can be
replaced by that on  $\R^n$. This, together with \cite{xu}, will
imply the validity of {\sf ii)} $\Rightarrow$ {\sf i)} in  both
Theorems \ref{lusin semigroups} and \ref{lusin type semigroups}
for the particular case of the Poisson semigroup on $\R^n$. This
is done by a careful analysis of the Poison kernels on $\T$ and on
$\R$ and a comparison of its essential parts.

\medskip

We shall also need a lemma, similar to Lemma \ref {restriction g
in the torus 1},  for the Poisson kernel on $\R$. We leave its
elementary proof to the reader.

\begin{lemma}\label{restriction g in R}
Let $\B, p, q$ and $\delta$ be as in the previous lemma. Then for
any $f\in L^p_\B(\R)$
$$
\Big\|\big[t\frac{\partial P_t}{\partial t}\ind_{(0,\delta)}(t)
\ind_{(\delta,\infty)}(|x|)\big]
*f\Big\|_{L^p_{L^q_\B((0,\infty),\frac{dt}{t})}(\R)}\le C_\delta
\|f\|_{L^p_\B(\R)}\ .
$$
\end{lemma}

\medskip

Now we state our result on the Lusin cotype for the Poisson
semigroup on $\R^n$.

\begin{theorem}\label{cotype in R}
Let $\B$ be a Banach space and $2\le q<\infty$. Then the following
statements are equivalent:
\begin{itemize}
\item[{\sf i)}] $\B$ is of Lusin cotype $q$.
\item[{\sf ii)}] For every (or equivalently, for some) positive integer
$n$ and for every (or equivalently, for some) $p\in(1,\infty)$
there is a constant $C>0$
 $$\|\g_q(f)\|_{L^p(\R^n)}\leq C\|f\|_{L^p_\B(\R^n)}\ ,\quad
 \forall\;  f\in L_\B^p(\R^n).$$
\item[{\sf iii)}] For every (or equivalently, for some) positive integer
$n$ there is a constant $C>0$ such that
$$\|\g_q(f)\|_{L^{1,\infty}(\R^n)}\leq C\|f\|_{L^1_\B(\R^n)}\ ,\quad
 \forall\;  f\in L_\B^1(\R^n).$$
\end{itemize}
The same equivalences hold with $\g^1_q$ or $\g^2_q$ instead of
$\g_q$ in {\sf ii)} and {\sf iii)}.
\end{theorem}

\noindent {\sc Proof}. In virtue of Proposition \ref{time-space},
we need only to prove the theorem for $\g^1_q$.

{\sf i)} $\Rightarrow$ {\sf ii)}. This is a particular case of
Theorem \ref{lusin semigroups}.

{\sf ii)} $\Leftrightarrow$ {\sf iii)}. This equivalence for a
given integer $n$ is already contained in Corollary \ref
{littlewood-paley}.

{\sf ii)} for $n>1$ $\Rightarrow$ {\sf ii)} for $n=1$. By
Corollary \ref{littlewood-paley}, it is enough to get the
boundedness from $L^\infty_{c,\B}$ into $\BMO$ of $\g^1_q$ on $\R$
from the same boundedness property of $\g^1_q$ on $\R^n$. To this
end, consider $\tilde x= (x_2,\dots, x_n)\in\R^{n-1}$, and $h\in
L^\infty_{c,\B}(\R)$, and define
$f(x)=h(x_1)\ind_{[0,1]^{n-1}}(\tilde x)$, where
$x=(x_1,x_2,\dots, x_n)\in\R^n$. The symmetric diffusion semigroup
generated by the Laplacian on $\R^n$ is given by convolution with
the Gaussian density. Then we have
\begin{eqnarray*}
 \Tt_tf(x)
 &=&\int_{\R^n}\frac{1}{(4\pi t)^{n/2}}
 e^{-\frac{|x-y|^2}{4t}}f(y)\,dy\\
 &=& C_0\int_{\R}\frac{1}{(4\pi t)^{1/2}}
 e^{-\frac{|x_1-y_1|^2}{4t}}h(y_1)\,dy_1 =C_0\Tt^1_th(x_1),
\end{eqnarray*}
where $\Tt^1_t$ is the heat kernel in $\R$. If we denote by
$\Pt^1_t$ the Poisson semigroup subordinated to $\Tt^1_t$ on $\R$
and by $P^1_t$ the Poisson kernel on $\R$, the formula
(\ref{subordinated semigroup}) implies that
$P_t*f(x)=\Pt_tf(x)=C_0\Pt^1_th(x_1)=C_0P^1_t*h(x_1)$, and
therefore $\g^1_qf(x)=C_0\g^1_qh(x_1)$. Now, for every interval
$I\subset\R$ consider $Q=I^n$ the cube in $\R^n$ whose sides are
the interval $I$. Then,
$$
\frac{1}{|Q|}\int_Q\g^1_qf(x)\,dx
=
\frac{1}{|I|^n}\int_{I^n}C_0\g^1_qh(x_1)\,dx_1\dots\,dx_n
=\frac{C_0}{|I|}\int_I\g^1_qh(x_1)\,dx_1.
$$
Therefore, and also by using similar arguments,
$$
\frac{1}{|Q|}\int_Q\big|\g^1_qf(x)-(\g^1_qf)_Q\big|\,dx
=\frac{C_0}{|I|}\int_I\big|\g^1_qh(x_1)-(\g^1_qh)_I\big|\,dx_1.
$$
Hence,
$$
\|\g^1_qh\|_{\BMO(\R)}=\frac1{C_0}\|\g^1_qf\|_{\BMO(\R^n)}\leq C\|f\|_{L^\infty_\B(\R^n)}
=C \|h\|_{L^\infty_\B(\R)}.
$$

{\sf ii)} for $n=1$ $\Rightarrow$ {\sf i)}. By Corollary
\ref{littlewood-paley} and the corresponding result in \cite{xu}
for the torus, we need to prove that for every $f\in L_\B^q(\T)$
$$\|G^1_q(f)\|_{L^q(\T)}\le C \|f\|_{L_\B^q(\T)}\ .$$
(Note that we take $p=q$ here; also recall that $G^1_q$ is the
$G_q$-function on the torus relative to the derivative in the
radius.) But by Lemma  \ref{restriction g in the torus 1} it is
enough to show that for   $\delta>0$ very close to 1
\begin{equation}\label{essential torus}
\Big\|\big[(1-r)\frac{\partial P_r(\varphi)}{\partial
r}\ind_{(1-\delta,1)}(r) \ind_{(-\delta,\delta)}(\varphi)\big]
*f\Big\|_{L^q_{L^q_\B((0,1),\frac{dr}{1-r})}(\T)} \leq C
\|f\|_{L^q_\B(\T)}.
\end{equation}
By the change of variables $r=e^{-t}$, we have
 $$\frac{\partial P_r(\varphi)}{\partial r}\ind_{(1-\delta,1)}(r)
 =\frac{2(1-e^{-t})^2-4(1+e^{-2t})\sin^2(\varphi/2)}
 {\big[(1-e^{-t})^2+4e^{-t}\sin^2(\varphi/2)\big]^2}\equiv k_t(\varphi).$$
Thus (\ref{essential torus}) is reduced to
\begin{equation}\label{essential torus1}
\int_\T \int_0^\varepsilon\big\|t\int_{-\delta}^{\delta}
 k_t(\varphi)
f(\theta-\varphi)\,d\varphi\big\|^q \,\frac{dt}{t} \, d\theta \leq
C^q \|f\|^q_{L^q_\B(\T)}\ .
\end{equation}
where $\varepsilon=\log\frac{1}{1-\delta}$. It is elementary to
decompose  $k_t(\varphi)$  as follows
 $$k_t(\varphi)=k^0_t(\varphi) +
 k^1_t(\varphi)+k^2_t(\varphi),$$
where
 \begin{equation}\label{k0}
 k^0_t(\varphi)=2\,\frac{t^2-\varphi^2}{(t^2+\varphi^2)^2}\,
 \ind_{(0,\varepsilon)}(t) \ind_{(-\delta,\delta)}(\varphi)
 \end{equation}
and where $k^1_t(\varphi)$ and $k^2_t(\varphi)$ are  supported on
$(0,\varepsilon)\times (-\delta, \delta)$ and satisfy
 $$|k^1_t(\varphi)|\le C_\delta \,\frac{t}{t^2+\varphi^2}\ ,
 \quad \ |k^2_t(\varphi)|\le C_\delta .$$
The verification of this decomposition, though entirely
elementary, could be tedious. One way to do this is to replace
each time only one term of $e^{-t}$ and $\sin(\varphi/2)$ by their
respective equivalents $1-t$ and $\varphi/2$ in $k_t(\varphi)$ and
in the functions so successively obtained. At each stage the
difference between the old and new functions is of type
$k^1_t(\varphi)$ when $e^{-t}$ is replaced or of type
$k^2_t(\varphi)$ when $\sin(\varphi/2)$ is replaced.

It is evident that
$$
\int_\T \int_0^\varepsilon\big\|t\int_{-\delta}^{\delta}
 k^2_t(\varphi)
f(\theta-\varphi)\,d\varphi\big\|^q \,\frac{dt}{t} \, d\theta \leq
C_\delta^q \|f\|^q_{L^q_\B(\T)}\ .
$$
It is also easy to get such an inequality for $k^1_t(\varphi)$.
Indeed, we have
\begin{eqnarray*}
 \int_\T \int_0^\varepsilon\big\|t\int_{-\delta}^{\delta}k^1_t(\varphi)
 f(\theta-\varphi)\,d\varphi\big\|^q \,\frac{dt}{t}\, d\theta
 &\leq& C_\delta^q \int_0^\varepsilon t^q\big[\int_{-\delta}^{\delta}
 \frac{t}{t^2+\varphi^2}\,\|f\|_{L^q_\B(\T)}\,d\varphi\big]^q\,
 \frac{dt}{t}\\
 &\le& C_{q,\delta}^q \|f\|^q_{L^q_\B(\T)}\ .
\end{eqnarray*}
Therefore, (\ref{essential torus1}) is reduced to
\begin{equation}\label{essential torus2}
\int_\T \int_0^\varepsilon\big\|t\int_{-\delta}^{\delta}
 k^0_t(\varphi)
f(\theta-\varphi)\,d\varphi\big\|^q \,\frac{dt}{t} \, d\theta \leq
C^q \|f\|^q_{L^q_\B(\T)}\ .
\end{equation}
Now we use the Poisson kernel $P_t$ on $\R$. By the definition of
$k^0_t$ in (\ref{k0}),
 $$k^0_t(x)=\frac{1}{2}\, \frac{\pa P_t(x)}{\pa t}
 \ind_{(0,\varepsilon)}(t)\ind_{(-\delta,\delta)}(x).$$
Put $\tilde f(x)=f(x)\ind_{(-\pi,\pi)}(x)$ for $x\in\R$. Then by
Lemma \ref{restriction g in R}, we see that (\ref{essential
torus2}) is further reduced to
 $$\int_\R \int_0^\varepsilon\big\|t
 \frac{\pa P_t}{\pa t}*\tilde f(x)\big\|^q \,\frac{dt}{t}\, dx
 \leq C^q \|\tilde  f\|^q_{L^q_\B(\R)}\ .$$
This last inequality follows from hypothesis {\sf iii)}.
Therefore, $\B$ is of Lusin cotype $q$, and thus the theorem is
proved. \cqd

\medskip

The following is the dual version of Theorem \ref{cotype in R}.

\begin{theorem}\label{type in R}
Let $\B$ be a Banach space and $1< q\le 2$. Then the following
statements are equivalent:
\begin{itemize}
\item[{\sf i)}] $\B$ is of Lusin type $q$.
\item[{\sf ii)}] For every (or equivalently, for some) $n\ge1$ and for every
(or equivalently, for some) $p\in(1,\infty)$ there is a constant
$C>0$
 $$\|f\|_{L^p_\B(\R^n)}\le C\,\|\g_q(f)\|_{L^p(\R^n)}\ ,\quad
 \forall\;  f\in L_\B^p(\R^n).$$
\end{itemize}
The same equivalence holds with $\g^1_q$ or $\g^2_q$ instead of
$\g_q$ in {\sf ii)}.
\end{theorem}

\noindent {\sc Proof}. {\sf i)} $\Rightarrow$ {\sf ii)} is a
particular case of Theorem \ref{lusin type semigroups}. {\sf ii)}
$\Rightarrow$ {\sf i)} is done by duality in virtue of Theorems
\ref{type-cotype} and \ref{type in R}. \cqd


\section{Ornstein-Uhlenbeck semigroup}

Our purpose of this section is to extend the results in the
previous one to the Poisson semigroup subordinated to the
Ornstein-Uhlenbeck semigroup on $\R^n$. Recall that this latter
semigroup is defined by
 $$O_tf(x)=\frac{1}{\big(\pi(1-e^{-2t})\big)^{n/2}}
 \int_{\R^n}\exp\big[-\,\frac{|e^{-t}x-y|}{1-e^{-2t}}\big]\;f(y)dy\
 .$$
We denote by $\big\{{\mathcal O}_t\big\}_{t\ge 0}$ the Poisson
semigroup subordinated to $\big\{{O}_t\big\}_{t\ge 0}$ as defined
in (\ref{subordinated semigroup}).

Let $\B$ be a Banach space and $1<q<\infty$. As for the usual
Poisson semigroup on $\R^n$, we introduce the Littlewood-Paley
$g$-function associated to $\big\{{\mathcal O}_t\big\}_{t\ge 0}\;$:
$$
g_q(f)(x)=\Big(\int_0^\infty t^q\|\nabla {\mathcal O}_t
f(x)\|_{\ell^2_\B}^q\,\frac{dt}{t}\Big)^{1/q},\quad x\in\R^n\ .
$$
Here $\nabla$ still denotes the gradient in  $\R^n\times \R_+$. We
shall also consider its two variants corresponding to the time
derivative and the space variable gradient, respectively:
$$
g^1_q(f)(x)=\Big(\int_0^\infty t^q\big\|\frac{\partial {\mathcal
O}_tf}{\partial t}(x)\big\|_{\B}^q\,\frac{dt}{t}\Big)^{1/q}\
$$
and
$$
g^2_q(f)(x)=\Big(\int_0^\infty t^q\|\nabla_x {\mathcal O}_t
f(x)\|_{\ell^2_\B}^q\,\frac{dt}{t}\Big)^{1/q} \ .
$$

The following is the analogue of Theorem \ref{cotype in R} for the
Ornstein-Uhlenbeck semigroup. $\gamma_n$ stands for the Gaussian
measure on $\R^n$, i.e., $\gamma_n=\exp(-|x|^2)dx$.

\begin{theorem}\label{cotype Orstein-Uhlenbeck}
Let $\B$ be a Banach space and $2\le q<\infty$. Then the following
statements are equivalent:
\begin{itemize}
\item[{\sf i)}] $\B$ is of Lusin cotype $q$.
\item[{\sf ii)}] For every (or equivalently, for some) positive integer
$n$ and for every (or equivalently, for some) $p\in(1,\infty)$
there is a constant $C>0$
 $$\|g_q(f)\|_{L^p(\R^n,\gamma_n)}\leq C\|f\|_{L^p_\B(\R^n,\gamma_n)}\ ,\quad
 \forall\;  f\in L_\B^p(\R^n,\gamma_n).$$
\item[{\sf iii)}] For every (or equivalently, for some) positive integer
$n$ there is a constant $C>0$ such that
$$\|g_q(f)\|_{L^{1,\infty}(\R^n,\gamma_n)}\leq C\|f\|_{L^1_\B(\R^n,\gamma_n)}
\ ,\quad \forall\;  f\in L_\B^1(\R^n,\gamma_n).$$
\end{itemize}
The same equivalences hold with $g^1_q$ or $g^2_q$ instead of
$g_q$ in {\sf ii)} and {\sf iii)}.
\end{theorem}

We also have a similar result for Lusin type.

\begin{theorem}\label{type Orstein-Uhlenbeck}
Let $\B$ be a Banach space and $1< q\le 2$. Then the following
statements are equivalent:
\begin{itemize}
\item[{\sf i)}] $\B$ is of Lusin type $q$.
\item[{\sf ii)}] For every (or equivalently, for some) $n\ge1$ and for every
(or equivalently, for some) $p\in(1,\infty)$ there is a constant
$C>0$
 $$\|f\|_{L^p_\B(\R^n, \gamma_n)}\le C\,\big(\big\|\int_{\R^n}f\,d\gamma_n\big\|_\B+
 \big\|g_q(f)\big\|_{L^p(\R^n, \gamma_n)}\big)\ ,\quad
 \forall\;  f\in L_\B^p(\R^n, \gamma_n\;).$$
\end{itemize}
The same equivalence holds with $g^1_q$ or $g^2_q$ instead of
$g_q$ in {\sf ii)}.
\end{theorem}

The proofs of the theorems above can be reduced to those on the
usual Poisson semigroup on $\R^n$ already considered in the
previous section. The usual technique dealing with operators related to
the Ornstein-Uhlenbeck semigroup consists in decomposing $\R^n$ into
two regions: one where the Gaussian and Lebesgue's measure are equivalent,
and the corresponding operators comparable, and the other where the kernels
of the operators can be estimated by a well behaved positive kernel.
This technique was invented by Muckenhoupt in the one-dimensional case, and
extended by Sj\"ogren to higher dimensions, for the maximal operator.
For vector-valued functions, the technique has been developed in \cite{btv}, see
also the references therein. Following this, for the $g$-function operator,
define the domains in $\R^n\times \R^n$:
$$D_1=\big\{(x, y): |x-y|<\frac{n(n+3)}{1+|x|+|y|}\big\}\
\mbox{and}\
D_2=\big\{(x, y): |x-y|<\frac{2n(n+3)}{1+|x|+|y|}\big\}\ .$$
Let $\varphi$ be a smooth function on $\R^n\times\R^n$ which is
supported on $D_2$, equal to $1$ on $D_1$ and satisfies
$$\|\nabla_x\varphi(x,y)\| + \|\nabla_y\varphi(x,y)\|
\leq C|x-y|^{-1}\ .$$
Let $T$ be  a Calder\'on-Zygmund singular integral operator on
$\R^n$ with kernel $k(x,y)$ as described at the beginning of
section~4 (and satisfying the conditions {\sf a)} and {\sf b)}
there). We decompose $T$ into its local and global parts
 $$T_{\rm glob}f(x)=\int k(x,y)[1-\varphi(x,y)]f(y)dy\quad
 \mbox{and}\quad T_{\rm loc}=T - T_{\rm glob}.$$
Now we can apply this decomposition to our favorite Littlewood-Paley
$g$-functions. We get the corresponding operators $g_{q,{\rm
loc}}$, $g_{q,{\rm glob}}\ $... for the subordinated Poisson
Ornstein-Uhlenbeck semigroup, and $\g_{q,{\rm loc}}$, $\g_{q,{\rm
glob}}\ $... for the usual Poisson semigroup.
The proofs of Theorems \ref{cotype Orstein-Uhlenbeck}
and \ref{type Orstein-Uhlenbeck} are sketchy since the estimates
needed are rather technical and can be
obtained in a parallel way as done in \cite{btv}.

\medskip

\noindent {\sc Proofs of Theorems \ref{cotype Orstein-Uhlenbeck}
and \ref{type Orstein-Uhlenbeck}}. We shall use the following
known facts from \cite{btv}
 \begin{itemize}
\item[{\sf a)}] $\displaystyle g_{q,{\rm glob}}f(x)\le
 \int_{\R^n}Q_1(x,y)\|f(y)\|_\B\, dy,\ $
where $Q_1$ is a nonnegative kernel supported on $D_1^c$ such that
the associated integral operator is of weak type (1, 1) and of
strong type $(p,p)$ for every $p\in(1,\infty)$ with respect to the
Gaussian measure;
\item[{\sf b)}] $\displaystyle \big|g_{q,{\rm loc}}f(x)-
 \g_{q,{\rm loc}}f(x)\big|\le\int_{\R^n}Q_2(x,y)\|f(y)\|_\B\,dy,\ $
where $Q_2$ is a nonnegative kernel supported on $D_2$ such that
 $$\sup_x\int_{\R^n}Q_2(x,y)dy<\infty\quad\mbox{and}\quad
 \sup_y\int_{\R^n}Q_2(x,y)dx<\infty\ .$$
Consequently, the integral operator associated to $Q_2$ is of
strong type $(p,p)$ for every $p\in[1,\infty)$ with respect to
both Lebesgue and Gaussian measures;
\item[{\sf c)}] similar statements hold for $g_q^1$ and $g_q^2$ in
place of $g_q$.
\end{itemize}
Then, using Theorem \ref{cotype in R}, we can show Theorem
\ref{cotype Orstein-Uhlenbeck} as in \cite {btv}. We omit the
details.

Theorem \ref{type Orstein-Uhlenbeck} is dual to Theorem
\ref{cotype Orstein-Uhlenbeck} in the case of $g_q^1$ , because of
the general Theorem \ref{type-cotype semigroup}. Similar duality
results hold for $g_q$ and $g_q^2$ too. Indeed, using the facts
above, we get a projection result (concerning $g_q$ and $g_q^2$)
for the subordinated Ornstein-Uhlenbeck Poisson semigroup similar
to Lemma \ref{projection}. Then we deduce the desired duality
result on $g_q$ and $g_q^2$. We leave again the details to the
interested reader.\cqd


\section{Almost sure finiteness}

We have seen in the previous sections (and also in \cite{xu}) that
the Lusin cotype property is equivalent to the boundedness of the
various generalized Littlewood-Paley $g$-functions on
$L^p$-spaces. The following result shows that this is still
equivalent to an apparently much weaker condition on the
$g$-functions

\begin{theorem}\label{finiteness ae th}
Given a Banach space $\B$ and $q\in[2,\infty)$, the following
statements are equivalent:
\begin{itemize}
\item[{\sf i)}] $\B$ is of Lusin cotype $q$.
\item[{\sf ii)}] For any $f\in L^1_\B(\T)$, $G^1_qf(z)<\infty$ for
almost every $z\in\T$.
\item[{\sf iii)}] For any $f\in L^1_\B(\R^n)$, $\g^1_qf(x)<\infty$ for
almost every $x\in\R^n$.
\end{itemize}
The equivalences also hold  when in statement {\sf ii)} $G^1_q$ is
replaced by $G^2_q$ or $G_q$, and also  in statement {\sf iii)}
$\g^1_q$ by $\g^2_q$ or $\g_q$.
\end{theorem}

\noindent{\sc Proof}. By  \cite{xu}, Theorem \ref{cotype in R} and
Corollary \ref{littlewood-paley}, we have  {\sf i)} $\Rightarrow$
{\sf ii)} and {\sf i)} $\Rightarrow$ {\sf iii)}. The two converse
implications are implicitly contained in \cite [VI.2] {ga-rdf:85}.
Let us first prove {\sf ii)} implies {\sf i)} (for $G^1_q$). To
this end, observe that
\begin{equation}\label{sup g torus}
G^1_q(f)(z)=\|Tf(z)\|_{L^q_{\B}((0,1),\frac{dr}{1-r})}
=\sup_{\varepsilon>0} \|T^\varepsilon
f(z)\|_{L^q_{\B}((0,1),\frac{dr}{1-r})}
\end{equation}
where $T^\varepsilon$ is the operator that sends $\B$-valued
functions to $L^q_{\B}((0,1),\frac{dr}{1-r})$-valued functions
given by
$$
T^\varepsilon f(z) =\big[(1-r)
\ind_{(\varepsilon,1-\varepsilon)}(r)\frac{\partial P_r}{\partial
r}\big]*f(z).
$$
It is clear that $T^\varepsilon$ is bounded from $L^1_\B(\T)$ to
$L^1_{L^q_\B((0,1),\frac{dr}{1-r})}(\T)$. Consequently, the
sublinear operator $f\mapsto
\|T^\varepsilon(f)\|_{L^q_\B((0,1),\frac{dr}{1-r})}$ is continuous
from $L^1_\B(\T)$ to $L^0(\T)$ (the latter space being equipped
with the measure topology).  By (\ref{sup g torus}), $G^1_q$ is
the supremum of these sublinear operators, and by {\sf ii)},
$G^1_q(f)\in L^0(\T)$ for all $f\in L^1_\B(\T)$. Therefore it
follows from the Banach-Steinhauss uniform continuity principle,
$G_q$ is continuous from $L^1_\B(\T)$ to $L^0(\T)$. Next, we apply
Stein's theorem. A proof of the scalar version can be found  in
\cite [VI.2] {ga-rdf:85}, and by using the ideas there, one can
prove the following vector-valued version.

\begin{lemma}\label{stein}
Let $G$ be a locally compact group with Haar measure $\mu$,  $\B$
be a Banach space  of Rademacher type $p_0$ and let
$T:L^p_\B(G)\longrightarrow L^0(G)$ be a continuous sublinear
operator invariant under left translations. Then  for every
compact subset $K$ of $G$  there exists a constant $C_K$ such that
$$
\mu(\{x\in K:\ |Tf(x)|>\lambda\})\leq
C_K\big(\frac{\|f\|_{L^p_\B}}{\lambda}\big)^{q}
$$
with $q=\inf\{p,p_0\}$. In particular, if the group $G$ is
compact, $T$ is of weak type $(p,q)$.
\end{lemma}

Let us recall that every Banach space is of Rademacher type $1$.
Then, $G^1_q$ is of weak type $(1,1)$, because it is clearly
sublinear and it is given by a convolution, which is invariant
under translations.

\medskip

The proof for the implication {\sf iii)} $\Rightarrow$ {\sf i)} is
similar. Again the sublinear operator $f\mapsto \g^1_q(f)$ is
continuous from  $L^1_\B(\R^n)$ to $L^0(\R^n)$. To infer as above
that it is of weak type (1,1), we use, instead of lemma
\ref{stein}, the following

\begin{lemma}\label{corollary of stein}
Let $\B$ be a Banach space of Rademacher type $p_0$ and $0<p\leq
p_0$. Then every translation and dilation invariant continuous
sublinear operator $ T:L^p_\B(\R^n)\longrightarrow L^0(\R^n)$ is
of weak type $(p,p)$.
\end{lemma}

This lemma can be proved in the same way as the corresponding
result in the scalar valued case in \cite [VI.2] {ga-rdf:85}. We
omit the details. Thus the proof of the theorem is  finished\cqd

\medskip

\noindent{\bf Remark.} Theorem \ref{finiteness ae th} holds also
for the $g$-function associated to the subordinated Poisson
Ornstein-Uhlenbeck semigroup.

\medskip

In the same spirit, we  also have a result similar to Theorem
\ref{finiteness ae th} in the case of martingales.

\begin{theorem}\label{martingale th}
Given a Banach space $\B$ and $2\le q<\infty$, the following
statements are equivalent:
\begin{itemize}
\item[{\sf i)}] $\B$ is of martingale cotype $q$.
\item[{\sf ii)}] If $f$ is a martingale bounded in $L^1_\B$, then
$S_q(f)<\infty$ almost everywhere.
\end{itemize}
\end{theorem}

For the proof of this theorem, we will use martingale transform
operators. Let $\B_1$ and $\B_2$ be two Banach spaces,
$(\Omega,\F,P)$ be a probability space, and $\{\F_n\}_{n\geq 1}$
be an increasing filtration of $\sigma$-subalgebras of $\F$. A
{\em multiplying sequence} $v=\{v_n\}_{n\geq 1}$ is a sequence of
random variables on $\Omega$ with values in ${\mathcal L}(\B_1,\B_2)$
such that each $v_n$ is $\F_{n-1}$-measurable and
${\displaystyle\sup_{n\geq 1}\norm{v_n}_{L^\infty_{{\mathcal
L}(\B_1,\B_2)}}<\infty}$. Given such a multiplying sequence,
define the martingale transform operator $T$ given by $v$ by
$(Tf)_n=\sum_{k=1}^n v_k\; d_kf$ for every martingale $f$. It is
proved in \cite{ma-to:00} that a martingale transform operator $T$
is of weak type $(1,1)$ iff it is of type $(p,p)$ for
$1<p<\infty$:
\begin{equation}\label{boundedness of T}
\sup_{\lambda>0}\lambda P\{(Tf)^\ast>\lambda\}\leq
C\,\norm{f}_{L^1_{\B_1}} \Leftrightarrow
\norm{(Tf)^\ast}_{L^p}\leq C_p\,\norm{f}_{L^p_{\B_1}},
\end{equation}
where it is understood that each side must hold for all martingales
$\{f_n\}$ with respect to the fixed filtration $\{{\mathcal F}_n\}$.
It is also proved there that, if $T$ is a translation invariant
martingale transform operator such that each term of its
multiplying sequence $\{v_k\}$ is a constant (in ${\mathcal
L}(\B_1,\B_2)$) and such that
\begin{equation}\label{reciproco2}
f^\ast\in L^1 \Longrightarrow \ Tf \mbox{ converges a.e.,}
\end{equation}
then $T$ also verifies the inequalities in (\ref{boundedness of
T}). By translation invariance of $T$ we mean that for any
$k_0\in\natt$, the sequence $\{v_k^{k_0}\}_{k\geq 1}$,
$v_k^{k_0}=v_{k_0+k}$, defines a martingale transform operator
 $T_{k_0}$ such that for any martingale $f$ bounded in $L^1_{\B_1}$,
 $$\norm{(Tf)_n}_{\B_2}=\big\|\sum_{k=1}^nv_kd_kf\big\|_{\B_2}
 =\big\|\sum_{k=1}^nv_{k_0+k}d_kf\big\|_{\B_2}
 =\big\|(T_{k_0}f)_n\big\|_{\B_2}\, ,\\ \forall\; n\ge 1.$$
Now, let  $Q_q$ be the martingale transform operator mapping
$\B$-valued martingales into $\ell^q_\B$-valued martingales
defined by the multiplying sequence $\{v_k\}_{k\geq 1}$  such that
each $v_k$ is the constant given by
$v_k(b)=(0,\stackrel{(k-1)}{\dots},0,b,0, \dots)$ for any $b\in
\B$. Then for a  $\B$-valued martingale $f$
 $$(Q_qf)_n=\sum_{k=1}^nv_kd_kf
 =(d_1f,d_2f,\dots,d_nf,0,\dots)\in\ell^q_\B$$
 $$\norm{(Q_qf)_n}_{\ell^q_\B}
 =\big(\sum_{k=1}^n\norm{d_kf}_\B^q\big)^{1/q}\ ,\quad
 (Q_qf)^\ast=\sup_n\norm{(Q_qf)_n}_{\ell^q_\B}=S_qf.$$

\noindent{\sc Proof of Theorem \ref{martingale th}}. {\sf i)}
$\Rightarrow$ {\sf ii)} is obvious. To prove the inverse, we use
that (\ref{reciproco2}) implies the inequalities (\ref{boundedness
of T}), applied to $Tf=Q_qf$.
 $Q_q$ is translation invariant and for $f^*\in L^1$,
$$
\|(Q_qf)_n-(Q_qf)_m\|_{\ell^q_\B}=
\Big(\sum_{k=m+1}^n\|d_kf\|_\B^q\Big)^{1/q} \longrightarrow
0\quad\mbox{a.e. as }n,m\rightarrow\infty,
$$
since it is the remaining of a convergent series (by {\sf
ii)}).\cqd

\medskip

We end with a final remark.

\medskip

\noindent{\bf Remark.} As in \cite {xu} for the torus, besides the
Littlewood-Paley $g$-function we can also consider the Lusin area
function on $\R^n$. In our vector-valued setting this function is
defined by
$$
\area_q(f)(x)=\Big(\int\int_{\Gamma(x)} t^q\|\nabla
P_t*f(y)\|_{\ell^2_\B}^q\,\frac{dydt}{t^{n+1}}\Big)^{1/q},
$$
where $\Gamma(x)$ is the cone with vertex $x$ and width $1$:
$$
\Gamma(x)=\{(t,y)\in\R^{n+1}_+:\ t> 0,\ |x-y|\leq t\}.
$$
Similarly, we can as well introduce the two variants $\area^1_q$
involving only the derivative in time and  $\area^2_q$ relative to
the gradient in the space variable. As in \cite {xu} for the
torus, all the preceding results (in sections 4 - 6) are still
valid with $\g_q$ replaced by $\area_q$. For instance, $\B$ is of
Lusin cotype $q$ iff $\|\area_q(f)\|_{L^p(\R^n)} \le C
\|f\|_{L_\B^p(\R^n)}$ for some (or all) $p\in(1,\infty)$, and iff
$f\in L_\B^1(\R^n) \Rightarrow \area_q(f)<\infty \;\mbox{a.e.
on}\; \R^n$.



\begin{thebibliography}{RdFRT}



\bibitem[BG]{bl-ga:87} Blasco, O., Garc{\'\i}a-Cuerva, J. $\ $
Hardy classes of Banach-space-valued distributions, {\em Math. Nachr.} {\bf 132}
(1987), 57--65.


\bibitem[BX]{bx} Blasco, O., Xu, Q., Interpolation between vector valued
Hardy spaces, {\em J. Funct. Anal.}, {\bf 102} (1991), 331--359.


\bibitem[Dav]{da:} Davies, E.B. {\em One-parameter diffusion semigroups}, L.M.S.
Monographs, {\bf 15}, Academic Press, London, 1980.







\bibitem[GR]{ga-rdf:85} Garc\'\i a-Cuerva, J., Rubio de Francia, J.L. $\ $
{\em Weighted norm inequalities and related topics}, North-Holland Mathematics
Studies, 116, Elsevier Science Publishers, 1985.

\bibitem[HTV]{btv} Harboure, E., Torrea, J.L., Viviani, B. $\ $
Vector valued extensions of operators related to the
Ornstein-Uhlenbeck semigroup, {\em J. Analyse Math.} (to appear).


\bibitem[Jou]{jou} Journ\'{e}, J.L $\ $ {\em Calder\'on-Zygmund operators,
Pseudo-differential operators and the Cauchy integral of Calder\'{o}n}, Lecture
Notes in Math. {\bf 994}, Springer-Verlag, Berlin, 1983.



\bibitem[LT]{LT} Lindenstrauss, J., Tzafriri, L. $\ $
{\em Classical Banach spaces II}, Springer-Verlag Berlin
Heidelberg New York, 1979.

\bibitem[MT]{ma-to:00} Martinez, T., Torrea, J.L., Operator-valued martingale
transforms {\em Tohoku Mathematical Journal}, {\bf 52},
(2000) 449--474.



\bibitem[Pi1]{pi:75} Pisier, G. $\ $ Martingales with values in uniformly convex
spaces, {\em Israel J. Math.}, {\bf 20}, (1975), 326-350.

\bibitem[Pi2]{pi:86} Pisier, G. $\ $  Probabilistic methods in the geometry of Banach
spaces, Lecture Notes in Mathematics, vol. {\bf 1206}, Springer-Verlag, Berlin.  1986.

\bibitem[Pi3]{pi:79} Pisier, G. $\ $ Some applications of the complex interpolation
method to Banach lattices, {\em J. Analyse Math.}, {\bf 35},
(1979), 264-281.

\bibitem[RdFRT]{rdfrt} Rubio de Francia, J.L., Ruiz, F.J., Torrea, J.L.  $\ $
Calder\'on-Zygmund theory  for operator valued kernels, {\em Adv. Math.},
{\bf 62} (1986), 7--48.


\bibitem[St1]{st:} Stein, E.M {\em Topics in Harmonic Analysis Related to the
Littlewood-Paley Theory}, Princeton University Press and the University of Tokio Press,
Princeton, New Jersey, 1970.

\bibitem[St2]{st:2} Stein, E.M {\em Singular Integrals and Differentiability Properties
of Functions}, Princeton University Press,
Princeton, New Jersey, 1970.



\bibitem[Xu]{xu} Xu, Q. $\ $ Littlewood-Paley theory for functions with values in
uniformly convex spaces, {\em J. reine angew. Math.} {\bf 504} (1998), 195--226.

\bibitem[Y]{yosido} Yosida, K {\em Functional Analysis},
Springer-Verlag, Berlin-Heidelberg-New York, 1968.


\end{thebibliography}
\end{document}